\documentclass{gtart}

\def\ifplaintex{\expandafter\ifx\csname documentclass\endcsname\relax}

\def\gtp{{\mathsurround=0pt\it $\cal G\mskip-2mu$eometry \&\ 
$\cal T\!\!$opology $\cal P\!$ublications}}  

\def\recd{{\small Received:\qua\receiveddate\ifx\reviseddate\relax
\else\qquad Revised:\qua\reviseddate\fi\par}} 

\def\lognumber#1{\def\thelognumber{#1}}
\def\volumenumber#1{\def\thevolumenumber{#1}}
\def\volumeyear#1{\def\thevolumeyear{#1}}
\def\papernumber#1{\def\thepapernumber{#1}}
\def\pagenumbers#1#2{\def\startpage{#1}\def\finishpage{#2}}
\def\published#1{\def\publishdate{#1}}

\def\received#1{\def\receiveddate{#1}}

\def\accepted#1{\def\accepteddate{#1}}
\def\asciititle#1{\def\theasciititle{#1}}

\long\def\asciiabstract#1{\long\def\theasciiabstract{#1}}
\def\asciikeywords#1{\def\theasciikeywords{#1}}

\let\\\par\let\thelognumber\relax\let\thevolumenumber\relax
\let\thepapernumber\relax\let\thevolumeyear\relax\let\startpage\relax
\let\finishpage\relax\let\publishdate\relax\let\receiveddate\relax
\let\reviseddate\relax\let\accepteddate\relax\let\theasciititle\relax
\let\theasciiauthors\relax
\let\theasciiabstract\relax\let\theasciikeywords\relax

\let\theasciiemail\relax

\ifplaintex
\font\logobig=cmssbx10 scaled 3836
\font\logomed=cmssbx10 scaled 2557
\else
\font\logobig=cmssbx10 scaled 4200
\font\logomed=cmssbx10 scaled 2800
\fi

\long\def\makeagttitle{   
\count0=\startpage
\agt\hfill      
\hbox to 45truept{\vbox to 0pt{\vglue -13truept{\logomed A\kern -.37em{\logobig 
T}\kern -.38em G}\vss}\hss}
\break
{\small Volume \thevolumenumber\ (\thevolumeyear)
\startpage--\finishpage\nl
Published: \publishdate}

\vglue .25truein

{\parskip=0pt\leftskip 0pt plus
1fil\def\\{\par\smallskip}{\Large\bf\thetitle}\par\medskip} \vglue
0.05truein

{\parskip=0pt\leftskip 0pt plus 1fil\def\\{\par}{\sc\theauthors}
\par\medskip}%
 
\vglue 0.03truein 

{\small\leftskip 25truept\rightskip 25truept{\bf Abstract}\stdspace\theabstract

{\bf AMS Classification}\stdspace\theprimaryclass
\ifx\thesecondaryclass\relax\else; \thesecondaryclass\fi\par
{\bf Keywords}\stdspace \thekeywords\par}\vglue 7truept

}   

\ifplaintex
\font\phead=cmsl9 scaled 950
\font\pnum=cmbx10 scaled 913
\font\pfoot=cmsl9 scaled 950
\headline{\vbox to 0pt{\vskip -4.5mm\line{\small\phead\ifnum
\count0=\startpage ISSN 1472-2739 (on-line) 1472-2747 (printed)
\hfill {\pnum\folio}\else\ifodd\count0\def\\{ }%
\ifx\theshorttitle\relax\thetitle\else\theshorttitle\fi\hfill{\pnum\folio}
\else\def\\{ and }{\pnum\folio}\hfill\ifx\theshortauthors\relax\theauthors
\else\theshortauthors\fi\fi\fi}\vss}}
\footline{\vbox to 0pt{\vglue 0mm\line{\small\pfoot\ifnum\count0=\startpage
\copyright\ \gtp\hfill\else
\agt, Volume \thevolumenumber\ (\thevolumeyear)\hfill\fi}\vss}}
\else
\font=cmsl9 scaled 1050
\font\lnum=cmbx10 
\font=cmsl9 scaled 1050
\makeatletter
\def\@oddhead{{\small\lhead\ifnum\count0=\startpage ISSN 1472-2739 
(on-line) 1472-2747 (printed)\hfill {\lnum\number\count0}\else\ifodd\count0
\def\\{ }\ifx\theshorttitle\relax \thetitle \else\theshorttitle\fi\hfill
{\lnum\number\count0}\else\def\\{ and }{\lnum\number\count0}
\hfill\ifx\theshortauthors\relax 
\theauthors\else\theshortauthors\fi\fi\fi}}\def\@evenhead{\@oddhead}
\def\@oddfoot{\small\lfoot\ifnum\count0=\startpage\copyright\ \gtp\hfill\else
\agt, Volume \thevolumenumber\ (\thevolumeyear)\hfill\fi}
\def\@evenfoot{\@oddfoot}
\makeatother
\fi
\let\maketitlepage\makeagttitle
\let\makeshorttitle\maketitlepage
\let\maketitle\maketitlepage

\newwrite\gtoutfile
\long\gdef\makeheadfile{  
{\def\\{, }\def\s{ }
\immediate\openout\gtoutfile head.xxx
\immediate\write\gtoutfile{To: math@arxiv.org}
\immediate\write\gtoutfile{Subject: put OR rep NNNNN:ppppp}
\immediate\write\gtoutfile{--text follows this line--}
\immediate\write\gtoutfile{Proxy-for: \ifx\theasciiauthors\relax
\theauthors\else\theasciiauthors\fi\s<\ifx\theasciiemail\relax\theemail\else\theasciiemail\fi>}
\immediate\write\gtoutfile{\noexpand\\}
\immediate\write\gtoutfile{Authors: \ifx\theasciiauthors\relax
\theauthors\else\theasciiauthors\fi}
{\def\\{ }\immediate\write\gtoutfile{Title: \ifx\theasciititle\relax
\thetitle\else\theasciititle\fi}}
\immediate\write\gtoutfile{Subj-class: GT or SG, GR etc}
\immediate\write\gtoutfile{MSC-class: \theprimaryclass\ifx\thesecondaryclass\relax\else, \thesecondaryclass\fi}
\immediate\write\gtoutfile{Journal-ref: Algebr. Geom. Topol. \thevolumenumber\s
(\thevolumeyear) \startpage-\finishpage}
\immediate\write\gtoutfile{Comments: Published by Algebraic and
Geometric Topology at}
\immediate\write\gtoutfile{\s\s\s  http://www.maths.warwick.ac.uk/agt/AGTVol\thevolumenumber/agt-\thevolumenumber-\thepapernumber.abs.html}
\immediate\write\gtoutfile{\noexpand\\}
\immediate\write\gtoutfile{}
\ifx\theasciiabstract\relax
\immediate\write\gtoutfile{\theabstract}\else
\immediate\write\gtoutfile{\theasciiabstract}\fi
\immediate\write\gtoutfile{}
\immediate\write\gtoutfile{\noexpand\\}
\immediate\write\gtoutfile{}
\immediate\closeout\gtoutfile}}  

\def\maketitlepage{\makeagttitle\makeheadfile}
\let\makeshorttitle\maketitlepage
\let\maketitle\maketitlepage

\lognumber{5}
\volumenumber{3}
\volumeyear{2003}
\papernumber{5}
\published{17 February 2003\newline Corrected: 20 January 2004 (see page 138)}
\pagenumbers{117}{145}
\received{16 November 2002}
\accepted{7 February 2003}

\usepackage{amsmath,amssymb,subfigure}
\usepackage[all,knot]{xy}
\usepackage{pstricks,pst-node,pstcol}

\newtheorem{thm}{Theorem}
 
 \newtheorem{cor}[thm]{Corollary}
 \newtheorem{lem}[thm]{Lemma}
 \newtheorem{prop}[thm]{Proposition}
 
 \theoremstyle{definition}
 \newtheorem{defn}[thm]{Definition}
 \theoremstyle{remark}
 \newtheorem{rem}[thm]{Remark}
 \numberwithin{equation}{section}

\newcommand{\CB}{\mathcal{B}}
\newcommand{\CC}{\mathcal{C}}
\newcommand{\CM}{\mathcal{M}}
\newcommand{\CN}{\mathcal{N}}
\newcommand{\CP}{\mathcal{P}}
\newcommand{\CS}{\mathcal{S}}

\newcommand{\FM}{\mathfrak{M}}
\newcommand{\Ga}{\alpha}
\newcommand{\Gb}{\beta}
\newcommand{\Gg}{\gamma}

\newcommand{\Gd}{\delta}

\newcommand{\Gk}{\kappa}
\newcommand{\Gl}{\lambda}

\newcommand{\Gr}{\rho}
\newcommand{\Gs}{\sigma}

\newcommand{\MB}{\mathbf}
\newcommand{\Au}{\operatorname{Aut}}
\newcommand{\n}{\newline}
\def\co{\colon\thinspace}

\begin{document}

\title{On $4$--fold covering moves}
\asciititle{On 4-fold covering moves}
\author{Nikos Apostolakis}
\address{Department of Mathematics, University of California\\Riverside 
CA 92521, USA}
\email{nea@math.ucr.edu}

\begin{abstract}
We prove the  existence of a finite set of moves sufficient
to relate any two representations of the same $3$--manifold 
as a $4$--fold simple branched covering of $S^3$. We also
prove a stabilization result: after adding a fifth trivial
sheet two local moves suffice. These results are 
analogous to results of Piergallini in degree $3$ and can be 
viewed as a second step in a program 
to establish similar results for arbitrary degree coverings
of $S^3$. 
    \end{abstract}

\asciiabstract{
We prove the  existence of a finite set of moves sufficient
to relate any two representations of the same 3-manifold 
as a 4-fold simple branched covering of S^3. We also
prove a stabilization result: after adding a fifth trivial
sheet two local moves suffice. These results are 
analogous to results of Piergallini in degree 3 and can be 
viewed as a second step in a program 
to establish similar results for arbitrary degree coverings
of S^3.}

\asciikeywords{Branched covering, covering move, colored braid, colored link, 
3-manifold}

\maketitle
\setcounter{section}{-1}

\section{Introduction}
\label{sec:intro}

It is known that every closed, orientable $3$--manifold is 
a simple branched covering of $S^3$ and furthermore the degree can be chosen
to be any number greater or equal to $3$. Such a covering can be represented 
by a ``colored link,'' that is, a link decorated with transpositions of some
symmetric group. Of course the same $3$--manifold admits many different 
representations as a colored link. 
In this paper we address the question of the  existence of a finite set of moves 
sufficient
to connect any two representations of the same $3$--manifold as a  colored link.
We concentrate on coverings of degree $4$, and give  two  answers.

In Theorem~\ref{thm:movesI} we give a set of seven moves 
sufficient to relate any two
presentations of the same manifold as a $4$-sheeted 
simple branched covering of the $3$-sphere. 

In Theorem~\ref{thm:main} we prove 
that two of those moves, the ``local moves,'' together with a form of
stabilization 
actually 
suffice.

These results and their proofs are modeled on similar results of Piergallini
for branched coverings of degree $3$
(\cite{Pier1991}, \cite{Pier1995}). It is unknown at this time if similar 
results are
true for coverings of degree greater than or
equal to $5$. Extensive computer calculations suggest that this is the 
case for 
degree $5$ and we conjecture that indeed Theorem~\ref{thm:main} is true for
arbitrary degree. If this conjecture is true we will have a combinatorial
description of $3$--manifolds comparable to the one given by Kirby calculus.
One application of such a calculus could be in 
the definition of new $3$-manifold invariants  or a geometric 
interpretation of existing  ones.

The theorems mentioned above appear in Section~\ref{sec:moves}.
The plan of the proofs  is the one introduced by Piergallini 
in ~\cite{Pier1991} and~\cite{Pier1995}:  Given $L_1$ and $L_2$, two colored
link representations  of the same $3$--manifold $M$, first  isotope the 
links so that they are in plat form. Each plat gives a  Heegaard splitting 
of $M$ whose
gluing homeomorphism is determined by lifting the  braid (viewed as 
a homeomorphism of the disk). A major part of our investigation therefore
 is to  
determine when two braids lift to the same gluing homeomorphism.
This is done in   
Section~\ref{sec:kernel}
 which culminates in  Theorem~\ref{thm:kernel}, result that 
which gives normal generators of the kernel of the lifting homomorphism.
The proof of this theorem is very similar to the proof of the corresponding theorem 
in~\cite{BirWajn1985}.

In the preceding Section~\ref{sec:complex}, the main technical section, we
determine the domain of the lifting homomorphism 
using the fact that it is the
stabilizer of the standard action of the braid group on the set of branched
coverings of the $2$-dimensional disk. The seemingly intractable calculations
are made possible by the following two observations. First,  the existence of 
 the exceptional
homomorphism $\xymatrix@1{{\CS_4}\ar[r]&{\CS_3}}$, which we call ``dimming the 
lights,'' means that
one  only needs to study the action of the $3$--fold stabilizer already computed 
in~\cite{BirWajn1985} instead of the action of the full braid group. Second, 
 the two ``local moves,'' can already be considered at this level as
moves between braids and therefore one only needs to compute in the ``reduced'' 
picture
i.e. after  taking quotient by the local moves.

Much  of the work in this paper was completed as part of the author's Ph.D.
thesis  under the supervision of Professor
Dennis Sullivan at the Graduate Center of CUNY. The author expresses  his
deep 
gratitude to  Dennis Sullivan for
his
guidance.      
The author would also like to thank the referee, whose careful reading and
insightful suggestions greatly improved the exposition of this paper.

\section{Preliminaries}
\label{sec:prel}
 We will always be working in the PL category. Thus in what follows manifold 
means PL manifold, 
homeomorphism means PL homeomorphism, submanifold means locally flat 
submanifold, etc.

We assume that the reader is familiar with the basic theory of branched 
coverings, especially in dimensions $2$ and $3$. An excellent
 background source 
  is~\cite{BernEdm1979} and the 
references contained therein. We also assume the reader is familiar 
with the contents of~\cite{BirWajn1985} as well  
as with  Piergallini's papers  regarding covering moves
 \cite{Pier1991}, \cite{Pier1995}. Finally we mention that the author's
thesis \cite{MyThesis} contains extensive background material. 
 
A $d$-fold {\em branched covering} is a map $p\co E \to B$
 between manifolds, that is a $d$--fold covering outside a co-dimension
 $2$ submanifold $L$ of $B$. $B$ is called the {\em base space},
$E$ the {\em total space}, $d$ is the degree and $L$ the 
{\em branching locus} of $p$.
In general, standard terminology of covering spaces (or more generally fiber 
bundles) will be used throughout the paper. 
Maps between branched coverings are defined as the obvious commutative
diagrams and we insist that an {\em isomorphism} covers the identity 
map of the base space.  
If $B$ has a basepoint~$*$, assumed to lie outside the branching locus,
 then a {\em labeled}  branched covering is a
 branched covering $p$ together with a {\em labeling}, i.e. 
 bijection of the fiber
$p^{-1}(*)$ with $\{1,\ldots,d\}$, where $d$ is the degree of $p$.
Maps between labeled  coverings are required to preserve the 
labellings. 

A covering of $B$ branched over $L$ is determined
(up to isomorphism) by the corresponding (unbranched) covering 
of $B \setminus L$. Therefore isomorphism classes of labeled
$d$-fold coverings  of $B$ branched over $L$ are in one one correspondence
with homomorphisms $\rho\co\pi_1(B\setminus L,*)\to\CS_d$ from the 
fundamental group of the complement of $L$ to the symmetric group on 
$d$-letters. This correspondence will be referred to as {\em monodromy}.
Quite often a covering will be confused with its monodromy; in particular
the same symbol will be used for both. 

If $p\co E\to B$ be a $d$-fold branched covering and $x \in E$,
 the {\em branching index} of $x$ is the ``number of sheets coming
together'' at $x$, or more formally the degree of $p$ restricted in 
a sufficiently small neighborhood of $x$. A branched covering 
will be called {\em simple} if the branching index of any point of 
$E$ is $1$ or $2$. When $B$ is simply connected 
the fundamental group $\pi_1(B\setminus L,*)$
is generated by {\em meridians},  simple closed paths that link
$L$ once. Then a branched covering is simple if the monodromy 
of each meridian is a transposition.
             
In the absence of  an explicit statement to the contrary,  ``covering''
means ``simple, labeled, branched covering''  throughout the 
article.  
    
\subsection{The Lifting Functor }
\label{sec:lift}
If $p\co E\to B$ is a covering branched over $L$ and $f\co X \to B$ is a
 homeomorphism then the pullback $f^*(p)$ is a covering of $X$ branched
 over the preimage $f^{-1}(L)$. Furthermore the equivalence class of
 $f^*(p)$ depends only on the isotopy class of $f$.
 In particular if 
$L$ is a fixed co-dimension $2$ submanifold of a manifold $B$, 
there is an action of the group of isotopy  classes of homeomorphisms for
$(B,L)$ on the set of isomorphism classes of coverings of 
$B$ branched over $L$. We are particularly interested in this 
action when the base space is the $2$-dimensional sphere $S^2$; it is 
technically more convenient however to consider coverings over
the $2$ dimensional disk $D^2$. 

In this case $L$ will consist of a finite set of points, say 
$L=\{A_0,..,A_{n-1}\}$ and a $d$--fold (simple) covering $\Gr$ 
branched over $L$
will be described by assigning a transposition $\Gs_i\in\CS_d$ to
each of the points $A_i$. Indeed the fundamental group of the 
punctured disk is 
freely 
generated by $\Ga_0,...,\Ga_{n-1}$, where each $\Ga_i$ is the concatenation 
of a path from from the basepoint to the boundary of 
a  small enough
disc containing  $A_i$ but no other branching points in its interior, the
boundary of that disc and the inverse of the path.
Given such a covering the monodromy $\Gs$ around the (positively oriented) 
boundary 
$\partial D^2$, will be equal to the product of all the monodromies 
$\Gs=\Gs_0\cdots\Gs_{n-1}$. Coverings over the sphere $S^2$  correspond
to coverings of $D^2$ with boundary monodromy $\Gs=\text{id}$, indeed the
total space of such a covering has $d$ boundary  components which 
can be ``filled in'' by $d$ disks mapping  with degree $1$ to the base $D^2$
resulting in a covering of $S^2$.    

It is well known that the mapping class group of 
$(D^2,L;\text{rel}\,\partial D^2)$ is isomorphic to $B_n$, the braid group
on $n$ strands. 
The isomorphism is given by assigning to the Artin generator 
$\Gb_i$ the counterclockwise rotation about an interval $x_i$ with endpoints 
$A_i$ and
$A_{i+1}$, for $i=0,\dotsc,n-2$. See~\cite{Bir1974} for details.
So there is a (right) action of the braid group on the set of coverings 
with $n$ 
branching values. 
The space of coverings branched over $L$ decomposes into orbits of this action 
according to the number of connected components of the total space
 and the boundary monodromy.
In particular this action is transitive on the set of coverings of $S^2$
with connected total space (see for example~\cite{BernEdm1979}).

\begin{figure}[htbp]\small
  \centering
  \begin{pspicture}(0,-.5)(10,2.5)
     \psline[linearc=.25](1,2)(1,1.5)(0,.5)(0,0)
     \psframe[linecolor=white,fillstyle=solid,fillcolor=white](.4,.9)(.6,1.1)
 \psline[linearc=.25](0,2)(0,1.5)(1,.5)(1,0) 
 \rput(0,2.3){$ij$}  \rput(0,-.3){$ij$} 
 \rput(1,2.3){$ij$}  \rput(1,-.3){$ij$} 
\rput(4,0){%
 \psline[linearc=.25](1,2)(1,1.5)(0,.5)(0,0)
     \psframe[linecolor=white,fillstyle=solid,fillcolor=white](.4,.9)(.6,1.1)
 \psline[linearc=.25](0,2)(0,1.5)(1,.5)(1,0)
 \rput(0,2.3){$ij$}  \rput(0,-.3){$ik$} 
 \rput(1,2.3){$jk$}  \rput(1,-.3){$ij$}
}%
\rput(8,0){%
 \psline[linearc=.25](1,2)(1,1.5)(0,.5)(0,0)
     \psframe[linecolor=white,fillstyle=solid,fillcolor=white](.4,.9)(.6,1.1)
 \psline[linearc=.25](0,2)(0,1.5)(1,.5)(1,0)
\rput(0,2.3){$ij$}  \rput(0,-.3){$kl$} 
 \rput(1,2.3){$kl$}  \rput(1,-.3){$ij$}
}%
   \end{pspicture}
  \caption{How a generator of $B_n$ acts}
  \label{fig:Baction}
\end{figure}
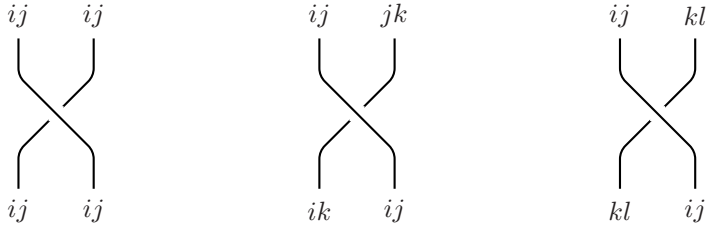

A combinatorial description of the braid group action can be given as
follows: Represent $\Gr$ as a coloring of $L$ by transpositions
as explained above.To see how a braid $\Gb$ acts on $\rho$ draw a diagram of
$\Gb$ with its top  endpoints coinciding with $L$, and let the
colors of $\rho$ ``flow down'' through the diagram according to
the rule that when a strand passes under another strand its color
gets conjugated by the color of the over strand. In this way we get a
new coloring of $L$ at the bottom of the diagram. This bottom
coloring represents $(\rho)\Gb$. Figure~\ref{fig:Baction} shows
the action of an Artin generator in the cases that the top colors coincide, 
``interact,'' or are
disjoint, respectively.

The theory of group actions is equivalent to the theory of 
groupoids. We find it convenient to think of the action of
the  braid group on sets of branched coverings in this framework. 
So given a natural number $d$, one can consider the groupoid 
of $d$--colored braids  whose objects are
colorings of $L$ by transpositions of the symmetric group $\CS_d$ and
its morphisms are colored braids. There is a {\em lifting functor} from
the groupoid of colored braids to the mapping class groupoid of
oriented surfaces. Restricting attention to a fixed
 object $\Gr$ one gets a lifting homomorphism  
$$\Gl\co\Au(\Gr)\to\FM(E(\Gr)),$$
from the group of automorphisms of $\Gr$ (i.e. in the language of
group actions the stabilizer of $\Gr$) to the mapping class
group of the total space of $\Gr$. Note that the elements of 
$\Au(\Gr)$ are called {\em liftable braids} in the literature.

Throughout this paper the Artin generators of the braid group $B_n$ will
be denoted by $\Gb_i$, for $i=0,\ldots,n-2$. Furthermore, for two braids
$\Gb$ and $\Gg$ we define $[\Gb]\Gg:=\Gg^{-1}\Gb\Gg$.  

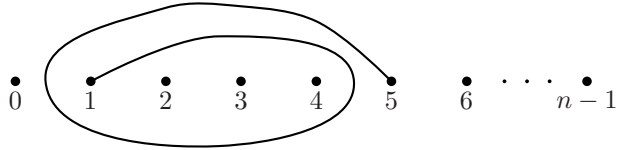
\begin{figure}[htp]\small
\centering
\begin{pspicture}(0,-.9)(10,1.1)
 \psdots[dotstyle=*,dotscale=1](0,0)(1,0)(2,0)(3,0)(4,0)(5,0)(6,0)(7.6,0)
 \pscurve (1,0)(2.7,.6)(4.5,0)(2.7,-.86)(.4,0)(2,.98)(3,1)(4,.8)(5,0)
 \psdots[dotstyle=*,dotscale=.4](6.5,0)(6.8,0)(7.1,0)
 \rput(0,-.26){0}
 \rput(1,-.26){1}
 \rput(2,-.26){2}
 \rput(3,-.26){3}
 \rput(4,-.26){4}
 \rput(5,-.26){5}
 \rput(6,-.26){6}
 \rput(7.6,-.26){$n-1$}
\end{pspicture}
\caption{The interval $d_4$}\label{fig:delta4}
\end{figure}

We will denote by  $\Gr^n(3)$  the ``standard''
 $3$-fold covering with $n$ branching points   
$$  (1,2), (1,2),  (2,3), \dotsc, (2,3)\quad,$$
 while its  automorphism group $\Au(\Gr^n)$ will be denoted by $L(n)$. 
 Birman and Wajnryb  studied $\Gr^n(3)$ in~\cite{BirWajn1985} (see 
also~\cite{BirWajn1994}). They found  that
$L(n)$, is generated by
$$\Gb_0, \Gb_1^3, \Gb_2,\Gb_3,\dotsc,\Gb_{n-2},\quad \text{and if } n\geq
6,\quad \delta_4$$
 where $\delta_4=[\Gb_4]\Gb_3\Gb_2\Gb_1^2\Gb_2\Gb_3^2\Gb_2\Gb_1$ is the
 rotation around the interval $d_4$ shown in Figure~\ref{fig:delta4}.

We will refer to the above 
generators as ``the Birman-Wajnryb generators'' or simply as ``the BW 
generators.''

\subsection{Dimming the lights}
\label{sec:dim}
\begin{figure}[htp]\small
\centering
\definecolor{lightblue}{rgb}{.5, .65, .9}
\definecolor{lightgreen}{rgb}{.5, .9, .7}
\definecolor{lightred}{rgb}{.9, .5, .5}
\begin{pspicture}(1,0)(3,3.2)
\rput(.85,0){1}
 \rput(2.5,3.2){3}
 \rput(4.15,0){2}
 \rput(2.57,1.1){4}
 \pnode(1,0){A}
 \pnode(2.5,3){B}
 \pnode(4,0){C}
 \pnode(2.75,.9){D}
 \ncline[linecolor=lightgreen]{A}{B}
 \ncline[linecolor=lightred]{A}{C}
 \ncline[linecolor=lightblue]{B}{C}
 \ncline[linestyle=dotted,linecolor=blue]{A}{D}
 \ncline[linestyle=dotted,linecolor=red]{B}{D}
 \ncline[linestyle=dotted,linecolor=green]{C}{D}
 \end{pspicture}
\caption{How to see the map $\kappa$}\label{fig:tetra}
\end{figure}
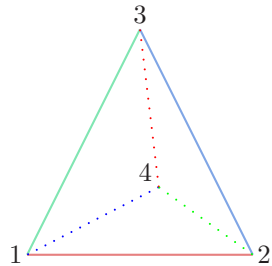

 There is an exact sequence of groups:
$$\xymatrix@1{ {1}\ar[r] & {V} \ar[r] & {\CS_4}
\ar[r]^{\kappa} & {\CS_3} \ar[r] & {1} }$$
 where $V=\{\text{id},(12)(34),(14)(23),(13)(24)\}$ is  Klein's
 four-group.
 Symmetric groups are generated by transpositions. Therefore, to describe $\Gk$  it 
suffices to describe its action on the  set of  transpositions of $\CS_4$.
Identify the edges of a numbered tetrahedron with the
transpositions of $\CS_4$ and the edges of its front face with the
transpositions of $\CS_3$ (see Figure~\ref{fig:tetra}). Then $\kappa$
fixes the front edges and sends each back edge to its opposite
edge.

 It is customary to identify, see~\cite{Pier1991} for example, the three
transpositions of $\CS_3$ with three colors Red, Green and Blue.
 The above exact sequence then
suggests the identification of the transpositions of $\CS_4$ with
three colors that come in light and dark shades. Denoting by
$\tilde X$ the dark shade of a color $X$ this identification is
given explicitly by:
$R=(12) $, $\tilde R=(34)$, $ B=(23)$, $ \tilde B=(14)$, $ G=(13)$
 and $\tilde G=(24)$.
The homomorphism $\kappa$ can then be described as ``dimming the
lights'' so that the two shades of the same color cannot be
distinguished.

 There is an induced functor from the groupoid of
``bi-tricolored'' braids to the groupoid of tricolored braids
which will also be referred to as dimming the lights. Dimming the
lights obviously commutes with the lifting functor.  In particular
if $\Gr$ is a $4$--fold covering then 
$$\Au(\Gr)\subset \Au\left(\Gk(\Gr)\right)$$
and furthermore the orbit of $\Gr$ under  $\Au\left(\Gk(\Gr)\right)$ 
is contained in $\Gk^{-1}\left(\Gk(\Gr)\right)$, the set of $4$--fold 
coverings that upon dimming the lights yield the same $3$--fold
 covering as $\Gr$.

\subsection{The reduced groupoid}
\label{sec:reduced}
\begin{figure}[htp]\small
\centering
\begin{pspicture}(-1,-4)(5,5)
\rput(0,2.5){%
\begin{xy}/r1.1cm/:
\vtwistneg\vcross\vtwistneg
\end{xy}
}%
\psline{<->}(1.2,2.5)(3,2.5)
\psline[linewidth=.4pt](3.5,.9)(3.5,4.1)
\psline[linewidth=.4pt](4.5,.9)(4.5,4.1)
\rput(0,-2){%
\begin{xy}/r1cm/:
\vtwistneg\vcross
\end{xy}
}%
\psline{<->}(1.2,-2)(3,-2)
\psline[linewidth=.4pt](3.5,-1)(3.5,-3)
\psline[linewidth=.4pt](4.5,-1)(4.5,-3)
\rput(2.1,2.9){$\mathcal{M}$}
\rput(-.6,4.4){$ij$} \rput(.5,4.4){$jk$}
\rput(3.5,4.4){$ij$} \rput(4.5,4.4){$jk$} 
\rput(-.9,3){$ik$} \rput(.8,3){$ij$}
\rput(-.9,2){$jk$} \rput(.8,2){$ik$}
\rput(-.6,.6){$ij$} \rput(.5,.6){$jk$}
\rput(3.5,.6){$ij$} \rput(4.5,.6){$jk$}
\rput(2.1,-1.6){$\mathcal{P}$}
\rput(-.6,-.7){$ij$} \rput(.5,-.7){$kl$}
\rput(3.5,-.7){$ij$} \rput(4.5,-.7){$kl$}
\rput(-.9,-2){$kl$}   \rput(.8,-2){$ij$}
\rput(-.6,-3.3){$ij$} \rput(.5,-3.3){$kl$}
\rput(3.5,-3.3){$ij$}  \rput(4.5,-3.3){$kl$} 
\end{pspicture}
\caption{The two local moves}\label{fig:localmoves}
\end{figure}
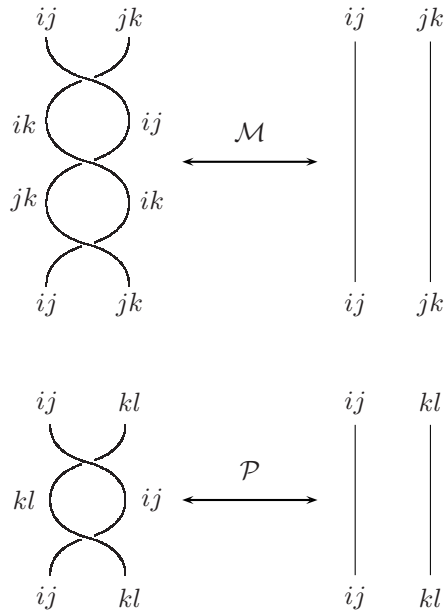

Consider the moves between colored braids shown in Figure~\ref{fig:localmoves}, 
where $i$, $j$, $k$, $l$ are distinct natural numbers.
 These moves are understood to  happen locally: The figure shows only 
the parts of the braids that are altered by the moves.    

 Let $\CN$ denote the equivalence relation
on (the set of morphisms of) the colored groupoid $\CB$ generated by the
two moves $\CM$ and $\CP$, that is two colored braids
$\Gb_1$ and $\Gb_2$ are equivalent modulo $\CN$ if and only if, they
can be related by a finite sequence of these moves. By the local
nature of the two moves it follows that $\CN$ is compatible with
the composition of colored braids and
thus there exists 
a quotient groupoid which we call the {\em reduced  groupoid} and denote 
by $\bar\CB$.
 Furthermore it is known that
\begin{lem}
\label{lem:liftred}
  The lifting functor factors through the reduced groupoid.  
\end{lem}
\begin{proof} It suffices to show that the braids on the 
left hand side of moves   $\CM$ and $\CP$ lift  to identity.
  Notice that the braid in the left hand side of move $\CM$ 
(respectively $\CP$)
is the third (respectively, second) power of a rotation around
an interval whose endpoints are labeled with interacting 
(respectively, disjoint) monodromies. Such an interval lifts
to a disjoint union of intervals and the third (respectively
second) power of the rotation around it lifts to a composition 
of the rotations around the lifted intervals. Such a composition
is isotopic to the identity. See~\cite{BernEdm1979} for details.    
\end{proof}
The induced functor on $\bar \CB$ will be denoted by $\bar{\Gl}$
and called the {\em reduced lifting functor} or, when no confusion
is likely, the lifting functor. 
If $\Gr$ is an object of the colored braid groupoid (i.e. a simple
covering) we denote by 
$N(\rho)$ the subgroup of $\Au(\rho)$ consisting of braids that are
$\CN$--equivalent to the identity morphism of $\Gr$. Notice that 
the local nature of the moves $\CM$ and $\CP$ implies that  $N(\rho)$
is a normal subgroup of  $\Au(\rho)$. The quotient group
   $$\frac{\Au(\Gr)}{N(\Gr)}$$
will be denoted by $\overline{\Au}(\Gr)$ and called the 
{\em reduced automorphism  group} of $\Gr$. Lemma~\ref{lem:liftred}
 specializes to the fact that there is a reduced lifting homomorphism
$$ \bar{\Gl}(\Gr) \co \overline{\Au}(\Gr) \to \FM\left( E(\Gr)\right).$$

\subsection{Coverings in dimension $3$}
\label{sec:3D}

In the case that the base space is the $3$--dimensional sphere $S^3$ the 
branching locus 
will be a {\em link} $L$. A $d$--fold covering branched over $L$ can be presented 
by assigning 
to each arc of a diagram of the link a transposition of $\CS_d$, in a way compatible 
with
 the Wirtinger
presentation. This means that at each crossing of the diagram exactly one 
of the three situations shown in Figure~\ref{fig:Baction} occurs. A link diagram
thusly decorated is called a $d$--colored link (diagram). One can easily trace
the effect of each of the Reidemeister moves on  a colored link diagram to
get a finite set of ``colored Reidemeister moves'' with the property that 
two colored link diagrams represent isomorphic coverings if and only if, they can be
related by a finite sequence of colored Reidemeister moves.    
The equivalence generated by colored Reidemeister moves will be called {\em colored
isotopy}. We say that the total space of the covering   is {\em represented} by 
the colored link. It is well known that every closed, oriented, $3$--manifold can 
be represented as a $d$--colored link for all $d\geq 3$ (see~\cite{Hil1974} 
or~\cite{Monte1974}).

\begin{figure}[htbp]
\centerline{\noindent\hglue 1.3cm
\mbox{\scalebox{.8}{
\begin{pspicture}(-.5,-.5)(15,5.5)
\psline[linearc=.15](.75,4)(.75,4.5)(1,5)(1.25,4.5)(1.25,1)(1.25,.5)(1,0)(.75,.5)(.75,4)
\psline[linearc=.15](2.75,4)(2.75,4.5)(3,5)(3.25,4.5)(3.25,1)(3.25,.5)(3,0)(2.75,.5)(2.75,4)
\psline[linearc=.15](3.75,4)(3.75,4.5)(4,5)(4.25,4.5)(4.25,1)(4.25,.5)(4,0)(3.75,.5)(3.75,4)
\psellipse(2.5,4.5)(2.4,.6)
\psellipse(2.5,.5)(2.4,.6)
\psframe[fillstyle=solid,fillcolor=white](.5,3.8)(4.5,1.2)
\rput(-.3,4.5){$B_{\text{up}}$}
\rput(-.3,.5){$B_{\text{low}}$}
\rput(-.3,2.5){$S^2\times I$}
\rput(2.5,2.5){\large{Braid}}
\psdots[dotsize=.05](1.8,4.5)(2,4.5)(2.2,4.5)(1.8,.5)(2,.5)(2.2,.5)
\psellipse(10.5,4.5)(2.4,.6)
\psellipse(10.5,.5)(2.4,.6)
\pscircle(9,4.5){.25} \pscircle(10.8,4.5){.25} \pscircle(11.8,4.5){.25}
\pscircle(9,.5){.25} \pscircle(10.8,.5){.25} \pscircle(11.8,.5){.25}
\psdots[dotsize=.05](9.7,4.5)(10,4.5)(10.3,4.5)(9.7,.5)(10,.5)(10.3,.5)
\psline{->}(10.5,3.5)(10.5,1.5)
\rput(11.3,2.5){$\lambda(\text{Braid})$}
\rput(13.3,4.5){$H_{\text{up}}$}
\rput(13.3,.5){$H_{\text{low}}$}
 \psline{->}(5.5,2.5)(7.5,2.5)
 \rput(6.5,2.8){Lift}
\end{pspicture}
}}}
\caption{How a colored  plat gives a  Heegaard splitting}\label{fig:Heeg}
\end{figure}
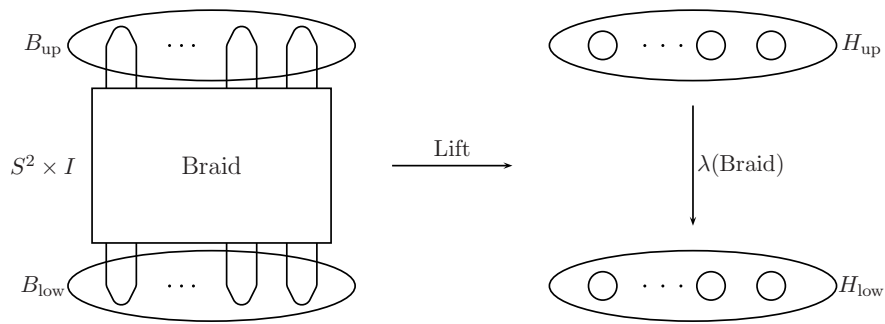

One way to understand the $3$--manifold represented by a colored link is via 
Heegaard splittings. Represent the link by a plat diagram and split $S^3$
as the union of a ball $B_{\text{up}}$ containing the caps of the plat, a tubular
neighborhood of its boundary containing the braid of the plat, and a second ball
$B_{\text{low}}$ containing the cups of the plat. Each ball lifts to a handlebody
and the braid lifts to a homeomorphism between the boundaries of the handlebodies,
and so we get a Heegaard splitting of the covering manifold. (See
Figure~\ref{fig:Heeg}, for more details consult  the references,
 for example~\cite{Pier1991}).

Moves $\CM$ and $\CP$ of Figure~\ref{fig:localmoves} (recall that  $i$, $j$,
 $k$, $l$ are distinct natural numbers)
 can be regarded as moves 
 on colored links. From now on they will collectively referred to as {\em the local
moves}. 
Since these two moves ``are in the kernel''  of the lifting
functor it is clear that two colored links that differ by a finite sequence of
moves represent the same $3$--manifold. There is also a {\em stabilization} move
which does not change the represented manifold but increases the degree of the
covering by $1$.  
If $L$ is a $d$--colored link, stabilization  adds an unknotted unlinked component 
colored
with a transposition of $\CS_{d+1}$ involving $d+1$.
That this move does not change the represented manifold is seen as
 follows. Let $M_d$ (respectively $M_{d+1}$) be the total space of the
 $d$--fold
(respectively $(d+1)$--fold) covering. In $M_{d+1}$ the preimage of $S^3$
 with a ball about the added component removed   consists of $M_d$
 with $d$ balls removed union $S^3$ with one ball removed. The preimage of
 the ball containing the added component consists of $d-1$ balls (filling in
 $d-1$ of the balls removed from $M_d$) union $S^2\times I$, joining the
 last puncture in $M_d$ to the punctured $S^3$, thus filling in the last
 puncture in $M_d$.   
 This move and its inverse 
will also be referred to
as an  \emph{ addition} or \emph{deletion of a trivial sheet}.

\section{The complex}
\label{sec:complex}

In this section we study the ``presentation complex'' associated 
with the action of $L(n)$, the automorphism group of the standard
$3$--fold covering $\Gr^n(3)$, on $\Gk^{-1}\left(\Gr^n(3)\right)$
 the set of $4$--fold
coverings that upon dimming the lights yield    $\Gr^n(3)$ (recall 
Section~\ref{sec:dim}). First we establish some notation:

\begin{defn}
  \label{def:notation}
For $\Gr\in \Gk^{-1}\left(\Gr^n(3)\right)$ we denote by $\CC(\Gr)$ the $L(n)$--orbit 
of $\Gr$.\newline
 The set of $4$-fold, {\em connected} coverings in $\Gk^{-1}(\Gr^n(3))$ 
that have their
first two monodromies equal to $(1,2)$ will be denoted by $\CC^n$, that is 
$$\CC^n:=\{\Gr\in\Gk^{-1}(\Gr^n(3)): 
\Gr(\Ga_0)=\Gr(\Ga_1)=(1,2)\quad \text{and $\Gr$ is surjective} \}.$$ 
Also for $\Gs\in\CS_n$ we denote by $\CC^n_\Gs$ the subset of $\CC^n$ 
consisting of 
coverings with boundary monodromy equal to $\Gs$, that is 
$$\CC^n_\Gs:=\{\Gr\in\CC^n : \mu(\Gr)=\Gs\},$$
where $\mu(\Gr)=\Gr(\Ga_0\Ga_1\dotsm\Ga_{n-1})$ denotes the boundary 
monodromy of
 $\Gr$. 
Finally, elements of $\CC^n $ are denoted according to which monodromies are
 equal to 
$(2,3)$, or alternatively which monodromies are equal to $(1,4)$:\newline
 for $I\subseteq
\{2,3,\dotsc,n-1\}$ define the $4$-fold coverings $\rho_I^n$,
$\tilde\rho_I^n$ of $D^2$ via:
$$\rho_I^n(\alpha_i):=\begin{cases} (12) & \text{if } i=0,1\\
                                                 (14) & \text{if } i\in I\\ 
                                                 (23) & \text{if } i\notin I
\end{cases} \hspace{1cm} 
\tilde\rho_I^n(\alpha_i):=\begin{cases} (12) & \text{if }i=0,1\\ 
                                       (23) & \text{if } i\in I\\ 
                                       (14) & \text{if } i\notin I
\end{cases}$$
\end{defn}

Next we look at the action of the  Birman-Wanjryb generators of $L(n)$.
The only generator whose action on $\CC^n$ is not obvious is $\Gd_4$.
 \begin{lem}
\label{lem:daction} 
$\Gd_4$ acts on $\CC^6$  as follows: it fixes $\Gr_{\emptyset}$,
 $\tilde\Gr_{\emptyset}$, $\Gr_{ij}$ and $\tilde\Gr_{ij}$, 
 while $(\Gr_i)\Gd_4=\tilde\Gr_i$.
\end{lem} 
\begin{proof} 
Since $\Gd_4\in L(6)$, it fixes   $\Gr_{\emptyset}$ and 
$\tilde\Gr_{\emptyset}$.  The 
remaining 
cases are easily checked by direct calculation. Two characteristic 
calculations
 are shown in 
Figure~\ref{fig:daction}.
\end{proof}
\begin{figure}[htp]
\centering
\mbox{\scalebox{.66}{
\begin{pspicture}(1,1)(5,10)
  \psline[linearc=.15](5,9.5)(5,6.3)(4,5)(1.1,4)(4.9,3)(1,2)(1,1.5)
 \psframe[linecolor=white,fillstyle=solid,fillcolor=white](4.3,5.47)(4.5,5.61)
 \psline[linearc=.15](1,9.5)(1,9)(4.9,8)(1.1,7)(3.6,6.3)(5,5)(5,1.5)
\psframe[linecolor=white,fillstyle=solid,fillcolor=white](1.9,8.65)(2.1,8.85)
\psframe[linecolor=white,fillstyle=solid,fillcolor=white](2.9,8.4)(3.1,8.6)
\psframe[linecolor=white,fillstyle=solid,fillcolor=white](3.9,8.15)(4.1,8.35)
\psframe[linecolor=white,fillstyle=solid,fillcolor=white](1.9,6.7)(2.1,6.9)
\psframe[linecolor=white,fillstyle=solid,fillcolor=white](2.9,6.4)(3.1,6.6)
\psframe[linecolor=white,fillstyle=solid,fillcolor=white](3.85,5.8)(4.15,6.05)
\psframe[linecolor=white,fillstyle=solid,fillcolor=white](1.9,4.2)(2.1,4.4)
\psframe[linecolor=white,fillstyle=solid,fillcolor=white](2.9,4.5)(3.1,4.7)
\psframe[linecolor=white,fillstyle=solid,fillcolor=white](3.9,4.9)(4.15,5.3)
\psframe[linecolor=white,fillstyle=solid,fillcolor=white](1.9,2.1)(2.1,2.4)
\psframe[linecolor=white,fillstyle=solid,fillcolor=white](2.9,2.4)(3.1,2.6)
\psframe[linecolor=white,fillstyle=solid,fillcolor=white](3.9,2.7)(4.1,2.9)
  \psline(2,9.5)(2,7.35)   \psline(2,7.12)(2,3.9)    \psline(2,3.65)(2,1.5) 
  \psline(3,9.5)(3,7.6)    \psline(3,7.4)(3,3.6)     \psline(3,3.4)(3,1.5)
  \psline(4,9.5)(4,7.9)    \psline(4,7.6)(4,3.35)    \psline(4,3.1)(4,1.5)
 \psdots(1,9.5)(2,9.5)(3,9.5)(4,9.5)(5,9.5)
        (1,1.5)(2,1.5)(3,1.5)(4,1.5)(5,1.5)
 \rput(1,9.8){$12$}
 \rput(2,9.8){$14$}
 \rput(3,9.8){$14$}
 \rput(4,9.8){$23$}
 \rput(5,9.8){$23$}
 \rput(1,1.2){$12$}
 \rput(2,1.2){$14$}
 \rput(3,1.2){$14$}
 \rput(4,1.2){$23$}
 \rput(5,1.2){$23$}
 \rput(1.1,7){$13$}
 \rput(1.1,4){$13$}
 \rput(3.4,5.1){$13$}
 \rput(3.5,2.38){$12$}
 \rput(2.5,8.9){$24$}
 \rput(2.5,6.89){$23$}
 \rput(4.3,7){$12$}
 \rput(1.7,5.8){$34$}
 \rput(2.7,5.8){$34$}
\end{pspicture}
\hspace{5cm}
\begin{pspicture}(1,1)(5,10)
  \psline[linearc=.15](5,9.5)(5,6.3)(4,5)(1.1,4)(4.9,3)(1,2)(1,1.5)
 \psframe[linecolor=white,fillstyle=solid,fillcolor=white](4.3,5.47)(4.5,5.61)
   \psline[linearc=.15](1,9.5)(1,9)(4.9,8)(1.1,7)(3.6,6.3)(5,5)(5,1.5)
\psframe[linecolor=white,fillstyle=solid,fillcolor=white](1.9,8.65)(2.1,8.85)
\psframe[linecolor=white,fillstyle=solid,fillcolor=white](2.9,8.4)(3.1,8.6)
\psframe[linecolor=white,fillstyle=solid,fillcolor=white](3.9,8.15)(4.1,8.35)
\psframe[linecolor=white,fillstyle=solid,fillcolor=white](1.9,6.7)(2.1,6.9)
\psframe[linecolor=white,fillstyle=solid,fillcolor=white](2.9,6.4)(3.1,6.6)
\psframe[linecolor=white,fillstyle=solid,fillcolor=white](3.85,5.8)(4.15,6.05)
\psframe[linecolor=white,fillstyle=solid,fillcolor=white](1.9,4.2)(2.1,4.4)
\psframe[linecolor=white,fillstyle=solid,fillcolor=white](2.9,4.5)(3.1,4.7)
\psframe[linecolor=white,fillstyle=solid,fillcolor=white](3.9,4.9)(4.15,5.3)
\psframe[linecolor=white,fillstyle=solid,fillcolor=white](1.9,2.1)(2.1,2.4)
\psframe[linecolor=white,fillstyle=solid,fillcolor=white](2.9,2.4)(3.1,2.6)
\psframe[linecolor=white,fillstyle=solid,fillcolor=white](3.9,2.7)(4.1,2.9)
  \psline(2,9.5)(2,7.35)   \psline(2,7.12)(2,3.9)    \psline(2,3.65)(2,1.5) 
  \psline(3,9.5)(3,7.6)    \psline(3,7.4)(3,3.6)     \psline(3,3.4)(3,1.5)
  \psline(4,9.5)(4,7.9)    \psline(4,7.6)(4,3.35)    \psline(4,3.1)(4,1.5)
 \psdots(1,9.5)(2,9.5)(3,9.5)(4,9.5)(5,9.5)
        (1,1.5)(2,1.5)(3,1.5)(4,1.5)(5,1.5)
 \rput(1,9.8){$12$}
 \rput(2,9.8){$23$}
 \rput(3,9.8){$23$}
 \rput(4,9.8){$23$}
 \rput(5,9.8){$14$}
 \rput(1,1.2){$12$}
 \rput(2,1.2){$14$}
 \rput(3,1.2){$14$}
 \rput(4,1.2){$14$}
 \rput(5,1.2){$23$}
 \rput(1.1,7){$13$}
 \rput(1.1,4){$24$}
 \rput(3.4,5.1){$24$}
 \rput(3.5,2.38){$12$}
 \rput(2.5,8.9){$13$}
 \rput(2.5,6.89){$23$}
 \rput(4.3,7){$12$}
 \rput(1.7,5.8){$12$}
 \rput(2.7,5.8){$12$}
\end{pspicture}
}}
\caption{How $\Gd_4$ acts}\label{fig:daction}
\end{figure}
\begin{cor}
 The  $L(n)$--action preserves $\CC^n$. 
\end{cor}
Now we can describe the orbits of the $L(n)$--action on $\CC^n$:  
\begin{prop}
  \label{prop:orbit}
For $\Gr\in\CC^n$, 
 $$\CC(\Gr)=\CC^n_{\mu(\Gr)} .$$ 
\end{prop}
\begin{proof} 
  Since the boundary monodromy is preserved by the braid action, 
 $\CC(\Gr)\subset\CC^n_{\mu(\Gr)}$. \newline
To prove the reverse inclusion notice that for  $\Gr\in\CC^n$, $\mu(\Gr)$ is equal to
\begin{itemize}
 \item id, if $n$ is even, and there is an even number of monodromies equal to $(1,4)$,
 \item $(1,4)(2,3)$, if $n$ is even and there is an odd number of monodromies equal to $(1,4)$,
 \item $(2,3)$,  if $n$ is odd and there is an even number of monodromies equal to $(1,4)$,
 \item $(1,4)$,  if $n$ is odd and there is an odd number of monodromies equal to $(1,4)$.
\end{itemize}
To see, for example, that every covering in $\CC_{\text{id}}^n$ 
 is in the $L(n)$-orbit of $\rho_{23}^n$ (where, of course, $n$ is even), observe
that the subgroup $\langle\Gb_2,\cdots,\Gb_{n-2}\rangle$ of $L(n)$, acts transitively on the
subset of $\CC_{\text{id}}^n$ consisting of
coverings with any fixed even number of monodromies equal to
$(14)$, while $\delta_4$, acting on suitably chosen coverings, increases
or decreases the number of monodromies equal to $(14)$ by two.\newline
The same argument works in the other cases as well.
\end{proof}

Next we enhance the orbits into $2$--complexes in a standard way:  
\begin{defn}
  For $\Gr\in \CC^n$ define a $2$-complex as follows
  \begin{itemize}
  \item The $0$-cells are the elements of $\CC(\Gr)$,
  \item the $1$-cells  are labeled by the Birman-Wajnryb 
 generators of $L(n)$. There is one $1$-cell labeled by $\Gb$ for each {\em unordered} pair $\{\rho',\rho''\}\subseteq\CC(\rho)$ with the property $\rho'\Gb=\rho''$,
   \item a path on the $1$--skeleton gives a word on the Birman-Wajnryb generators of $L(n)$.
         There is a $2$-cell  attached along each closed path whose word lies in $N(\Gr')$ for 
         some $\Gr'$ in the   path.    
  \end{itemize}
The resulting $2$--complex will be denoted by $\MB{C}(\rho)$. In general we
will use bold font to signify the transition from $L(n)$--sets to the
corresponding complexes.    
\end{defn}

Note that since all BW generators act as involutions we don't need to consider
 oriented edges. Clearly the fundamental group $\pi_1\left(\MB{C}(\rho),\rho\right)$
 is isomorphic to $\overline{\Au}(\rho)$. 

We are going to use the following terminology: a {\em vertex} 
is a $0$--cell, an {\em edge} is a $1$--cell attached to two 
{\em distinct}  vertices and a {\em loop} is a $1$--cell attached 
to a single vertex. A {\em simple} path is a path without 
 closed subpaths. 
A  {\em lasso} is a closed path of 
the form $wlw^{-1}$ where $l$ is a loop and $w$ is a simple
path, $w$ is the {\em tail} of the lasso and $l$ is its {\em head}. 

 Conjugation by $(12)(34)$ (i.e interchanging the dark and light
shades of Blue) defines an involution $\,\widetilde{}\,$ on
$\CC^n$. Clearly $\,\widetilde{}\,$ interchanges
$\CC_{(23)}^n$ and $\CC_{(14)}^n$ and preserves
$\CC_{\text{id}}^n$ and $\CC_{(14)(23)}^n$. In particular
$\MB{C}_{(23)}^n$ and $\MB{C}_{(14)}^n$ are homeomorphic. In what follows we will refer to
 $\,\widetilde{}\,$ simply as ``symmetry.''

 Given a covering with $n$ branch values one can ``add''
one more point with prescribed monodromy at the end to get a
covering with $n+1$ branch values. More formally, for each $n$
there are two embeddings:
$$ i_{(23)}^n\co\CC^n\hookrightarrow  \CC^{n+1},$$
and
$$i_{(14)}^n:\CC^n\hookrightarrow  \CC^{n+1}$$
which are obviously $L(n)$-equivariant if $B_n$ is (as usual)
considered to be the subgroup of $B_{n+1}$ that fixes the last
point. Therefore these embeddings extend to embeddings of the corresponding
complexes which we will still denote by $ i_{(23)}^n$ and $i_{(14)}^n$.

 \begin{prop}\label{prop:dec}
 For even $n$:
 \begin{itemize}
 \item[\rm(a)]At the level of $0$--skeletons:
$$\CC_{\text{id}}^n=i_{(23)}^{n-1}(\CC_{(23)}^{n-1})\bigsqcup 
                                 i_{(14)}^{n-1}(\CC_{(14)}^{n-1}),$$
and each vertex of $i_{(23)}i_{(14)}(\MB{C}_{(14)(23)}^{n-2})$ is
joined by an edge labeled by $\Gb_{n-2}$ to the corresponding
vertex of $i_{(14)}i_{(23)}(\MB{C}_{(14)(23)}^{n-2})$. These are the
only new edges of $\MB{C}_{\text{id}}^n$ .
\item[\rm(b)]  At the level of $0$--skeletons:
 $$\CC_{(14)(23)}^n=i_{(23)}^{n-1}(\CC_{(14)}^{n-1})\bigsqcup
 i_{(14)}^{n-1}(\CC_{(23)}^{n-1})$$
and each vertex of $i_{(23)}i_{(14)}(\MB{C}_{\text{id}}^{n-2}
\sqcup\{\rho_\emptyset^{n-2}\}\sqcup\{\tilde\rho_\emptyset^{n-2}\})$
is joined by an edge labeled by $\Gb_{n-2}$ to the corresponding
vertex of $i_{(14)}i_{(23)}(\MB{C}_{\text{id}}^{n-2}
\sqcup\{\rho_\emptyset^{n-2}\}\sqcup\{\tilde\rho_\emptyset^{n-2}\})$.
These are the only new edges of $\MB{C}_{(14)(23)}^n.$
\end{itemize}
 For odd $n$:
\begin{itemize}
\item[] At the level of $0$--skeletons
 $$\CC_{(23)}^n=i_{(23)}^{n-1}(\CC_{\text{id}}^{n-1}\sqcup \{\tilde\rho_\emptyset^{n-1}\})
 \bigsqcup i_{(14)}^{n-1}(\CC_{(14)(23)}^{n-1}),$$
and each vertex of $i_{(23)}i_{(14)}(\MB{C}_{(14)}^{n-2}
\sqcup\{\tilde\rho_\emptyset^{n-2}\})$ is joined by an edge
labeled by $\Gb_{n-2}$ to the corresponding vertex of
$i_{(14)}i_{(23)}(\MB{C}_{\text{id}}^{n-2}
\sqcup\{\tilde\rho_\emptyset^{n-2}\})$. These are the only new
edges of $\MB{C}_{(23)}^n$.
\end{itemize}
\end{prop}
\begin{proof}
  The proof is straightforward: remove the last point and examine 
 the result. 
\end{proof}

The complexes $\MB{C}(\Gr)$  have the property that every closed 
simple path is made up from squares.

\begin{prop} Every edge of $\MB{C}(\Gr)$ that is not part of  an isolated chain,
that is a path of the form $\xymatrix@1{{\Gr_1}\ar@{-}[r]&{\Gr_2}\ar@{-}[r]
&{\cdots}\ar@{-}[r]&{\Gr_m} }$ where all non displayed $1$--cells are loops,
is part of  a square
$$
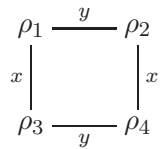
\xymatrix{ {\rho_1}\ar@{-}[d]_{x}\ar@{-}[r]^{y} &
{\rho_2}\ar@{-}[d]^{x}\\
{\rho_3}\ar@{-}[r]_{y} & {\rho_4} } $$
  Every such square bounds a $2$--cell.
\end{prop} 
\begin{proof} The first statement follows from the fact that 
all BW generators act as involutions. For the second statement, observe 
that for such a square either $x$ and $y$ are non-adjacent $\Gb$'s
or one of them is $\Gd_4$ and the other is different from $\Gb_5$.
In all these cases, except when $x=\Gd_4$ and $y=\Gb_4$ (or vice versa),  
 $x$ and $y$ commute in the braid group. In the remaining case, 
 according to Lemma~\ref{lem:comm}, $x$ and $y$ commute modulo the
subgroup $N$.  
\end{proof}

\begin{lem}\label{lem:comm} 
$\Gd_4$ commutes with $\Gb_4$ modulo the move  $\CM$.
\end{lem}
\begin{proof}
  This is proven (with a different formulation) in~\cite{BirWajn1985},
in the middle of page~37. For a pictorial proof consult~\cite{MyThesis}.
\end{proof}

 \subsection{Generators}
\label{sec:gen}
In this section we derive generators for the reduced automorphism
 group of  
the standard $4$--fold covering $\Gr^n_{23}$. Namely, we prove
 the following main technical result:

\begin{thm}
  \label{thm:gen}
 For all even $n\geq 6$, $\overline{\Au}(\rho_{23}^n)$ is
generated by $\Gb_0$, $\Gb_2$, $\Gb_4$, $\Gb_5$,
$\dotsc, \Gb_{n-2}$,
$\delta_4$ and, if $n\geq 8$, $\delta_6$, where,
$$\delta_6=[\delta_4]\Gb_5^{-1}\Gb_4^{-1}\Gb_3^{-1}\Gb_2^{-1}
\Gb_6^{-1}\Gb_5^{-1}\Gb_4^{-1}\Gb_3^{-1}.$$
\end{thm}

\begin{figure}[htp]
\centering 
\mbox{\scalebox{0.9}{\begin{pspicture}(0,-.9)(10,1.1)
 \psdots[dotstyle=*,dotscale=1](0,0)(1,0)(2,0)(3,0)(4,0)(5,0)(6,0)
                               (7,0)(8,0)(9.8,0) 
 \pscurve
 (1,0)(2,.6)(3.5,.8)(4.5,.6)(6.3,0)(4.6,-.4)(3.6,0)(2.5,.4)(1.5,0)
(1,-.4)(.7,0)(2,.8)(4,1)(6,.8)(7,0)
 \psdots[dotstyle=*,dotscale=.4](8.6,0)(8.9,0)(9.2,0)
  \rput(0,-.8){0}
 \rput(1,-.8){1}
 \rput(2,-.8){2}
 \rput(3,-.8){3}
 \rput(4,-.8){$4$}
 \rput(5,-.8){$5$}
 \rput(6,-.8){$6$}
 \rput(7,-.8){$7$}
 \rput(8,-.8){$8$}
 \rput(9.8,-.8){$n-1$}
\end{pspicture}}}
\caption{The interval $d_6$}\label{fig:delta6}
\end{figure}
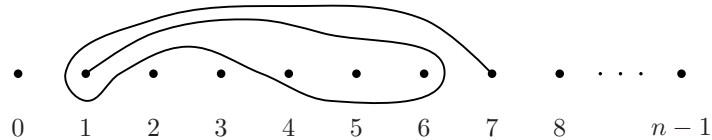

Note that $\Gd_6$ is rotation around the interval $d_6$ shown in 
figure~\ref{fig:delta6}.
Actually the following stronger statement will be proved:

\begin{thm}\label{thm:tech}
\begin{itemize}\item[\rm(a)] For each $n\geq 6$, 
$\overline{\Au}(\rho_{23}^n)$ is
generated by $\Gb_0$, $\Gb_2$, $\Gb_4$, $\Gb_5$,
$\dotsc, \Gb_{n-2}$,
$\delta_4$ and, if $n\geq 8$, $\delta_6$.
\item[\rm(b)] For each even $n\geq 6$, 
$\overline{\Au}(\rho_{2,3,n-1}^n)$ is
generated by $\Gb_0$, $\Gb_2$, $\Gb_4$,$\dots, \Gb_{n-3}$,
$[\Gb_{n-2}]\Gb_{n-3}^{-1}\dotsm\Gb_3^{-1}$ and if $n\geq8$,
$\delta_4$ and, if $n\geq 10$, $\delta_6$.

\end{itemize}
\end{thm}
  
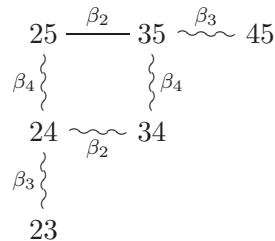
\begin{figure}[htp]
$$
\xymatrix{
 {25}\ar@{~}[d]_{\beta_4}\ar@{-}[r]^{\beta_2}
 &{35}\ar@{~}[d]^{\beta_4}\ar@{~}[r]^{\beta_3}
 &{45}\\
 {24}\ar@{~}[d]_{\beta_3}\ar@{~}[r]_{\beta_2}
 &{34}\\
 {23}         }
$$
\caption{$\CC(\rho_{23}^6)$ with a maximal tree}\label{fig:six1}
\end{figure}

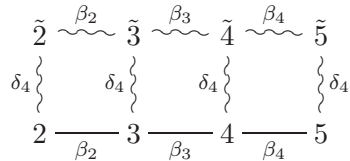
\begin{figure}[htp]
$$
\xymatrix{
 {\tilde 2}\ar@{~}[d]_{\delta_4}\ar@{~}[r]^{\beta_2}
 &{\tilde 3}\ar@{~}[d]_{\delta_4}\ar@{~}[r]^{\beta_3}
 &{\tilde 4}\ar@{~}[d]_{\delta_4}\ar@{~}[r]^{\beta_4}
 &{\tilde 5}\ar@{~}[d]^{\delta_4}\\
 {2}\ar@{-}[r]_{\beta_2}
 &{3}\ar@{-}[r]_{\beta_3}
 &{4}\ar@{-}[r]_{\beta_4}
 &{5}            }
$$
\caption{$\CC(\tilde\rho_4^{6})$ with a maximal
tree}\label{fig:six2}
\end{figure}

We start by establishing the following 
 \begin{lem}
   \label{lem:lasso} Let $\Gs$ be any of id, $(2,3)$, $(1,4)(2,3)$.  
The fundamental group of $\MB{C}_\Gs^n$ is generated by
 lassos. Furthermore the tail of a lasso is not important: two
lassos based at the same vertex and having the same head (at the
same vertex labeled by the same element) are equal.
 \end{lem}
\begin{proof}
Attach $2$--cells to kill all loops in $\MB{C}_\Gs^n$. It suffices to prove that  the
 resulting complex
 is simply connected. Proceed by induction on the number of branch
points $n$.

For $n=6$ a simple inspection of $\CC^6_{\text{id}}$ and $\CC^6_{(1,4)(2,3)}$ shown 
in 
Figures~\ref{fig:six1} and~\ref{fig:six2} confirms the truth of the claim. 
(Recall  Lemma~\ref{lem:comm}.)\n 
 For the induction step, notice that by Proposition~\ref{prop:dec}, 
$\MB{C}^n_{\Gs}$ is the union  
of simply connected complexes with connected intersection. Thus the claim follows 
from Van
Kampen's theorem. 
\end{proof}

\begin{proof}[Proof of Theorem~\ref{thm:tech}] 
Let $G$ be the subgroup of $\overline{\Au}(\Gr)$ generated by the 
listed elements. By Lemma~\ref{lem:lasso} it suffices to prove that
all lassos belong to $G$. 
Proceed by induction on the number of branch points $n$. For $n=6$, two simple applications 
of the 
Reidemeister-Schreier method using the transversals shown in
 Figures~\ref{fig:six1} and~\ref{fig:six2} prove the claim.\n
Assume then that the result has already been proved 
for all smaller values of~$n$.

If $n$ is odd,
\begin{align*} \CC(\rho_{23}^n) &=\CC_{(23)}^n\\
                                &=i_{(23)}(\CC_{\text{id}}^{n-1}\sqcup
                                 \{\tilde\rho_\emptyset^{n-1}\})
                              \bigsqcup i_{(14)}(\CC_{(14)(23)}^{n-1})
\end{align*}
and the lassos fall into three cases:
\begin{description}
\item[Case a] Lassos with heads at
$i_{(23)}(\MB{C}_{\text{id}}^{n-1})$
\end{description}
In this case, by the induction hypothesis, one needs only to 
check lassos with head labeled by $\Gb_{n-2}$.
  After writing
$$\CC_{\text{id}}^{n-1}=i_{(23)}(\CC_{(23)}^{n-2})
 \bigsqcup i_{(14)}(\CC_{(14)}^{n-2})$$
 one  observes that these lassos actually reside in
 $i_{(23)}i_{(23)}(\MB{C}_{(23)}^{n-2})$
and thus their  
 tail can be chosen
 to commute with $\Gb_{n-2}$. Therefore, all of these lassos
 are equal to $\Gb_{n-2}$.
\begin{description}
\item[Case b] Lassos with heads at
$i_{(14)}(\MB{C}_{(14)(23)}^{n-1})$
\end{description}
By the induction hypothesis one needs only to check the lassos
with heads at $\rho_{2,3,n-2,n-1}$, the lasso with head at
$\rho_{2,n-3,n-2,n-1}$ labeled by $\Gb_{n-3}$, and the lasso with
head at $\rho_{5,6,n-2,n-1}$ labeled by
$\delta_4$.

For the lassos with head at $\rho_{2,3,n-2,n-1}$ labeled by  $\Gb_i$ with
$i=2,4,5,\dotsc,n-5$ 
 or by $\delta_4$ one can use the
path $\Gb_{n-2}\Gb_{n-3}$ to pull them at $\rho_{2,3,n-3,n-2}\in
i_{(23)}(\MB{C}_{\text{id}}^{n-1})$ with the same label and refer to
case~a).
For the lasso with head at $\rho_{2,3,n-2,n-1}$ labeled by
 $\Gb_{n-2}$ use the path $\Gb_{n-2}\Gb_{n-3}$ to pull it at
 $\rho_{2,3,n-3,n-2}\in i_{(23)}(\MB{C}_{\text{id}}^{n-1})$ with label
 $\Gb_{n-3}$ and refer to case~a).

For the lasso with head at $\rho_{2,3,n-2,n-1}$ labeled by
 $\Gb_{n-4}$ one can use the path 
$\Gb_{n-3}\Gb_{n-2}\Gb_{n-4}\Gb_{n-3}$ to pull it at
 $\rho_{2,3,n-4,n-3}\in i_{(23)}(\MB{C}_{\text{id}}^{n-1})$ with label
 $\Gb_{n-2}$ and refer to case~a).

 For the lasso
with head at $\rho_{2,n-3,n-2,n-1}$ labeled by $\Gb_{n-3}$, one can use
the path $\Gb_{n-2}\Gb_{n-3}\Gb_{n-4}$ to pull it at
 $\rho_{2,n-4,n-3,n-2}\in i_{(23)}(\MB{C}_{\text{id}}^{n-1})$ with label
 $\Gb_{n-4}$ and refer to case~a).

 Finally for  the lasso with head at $\rho_{5,6,n-2,n-1}$ labeled by
$\delta_4$, use the path $\Gb_{n-2}\Gb_{n-3}$ to pull it at
 $\rho_{5,6,n-3,n-2}\in i_{(23)}(\MB{C}_{\text{id}}^{n-1})$ with label
 $\delta_4$ and refer to case~a).

\begin{description}
\item[Case c] Lassos with heads at
$\tilde\rho_{n-1}^n$
\end{description}
For those with heads labeled by $\delta_4$ or $\Gb_i$ with
$i=1,2,\ldots,n-4$ use $\Gb_{n-2}$ to pull them at
$\tilde\rho_{n-2}^n\in i_{(14)}(\MB{C}_{(14)(23)}^{n-1})$ with the
same
label, and refer to case~b).\n
For the one with head labeled by $\Gb_{n-3}$ use
$\Gb_{n-3}\Gb_{n-2}$ to pull it at $\tilde\rho_{n-3}^n\in
i_{(14)}(\MB{C}_{(14)(23)}^{n-1})$ with label $\Gb_{n-2}$, and refer
to case~b).\n
This concludes the proof in the odd $n$ case.

If  $n$ is even we have to check two different 
boundary monodromies. If the boundary monodromy is equal 
to identity then  
 \begin{align*}
\CC(\rho_{23}^n) & = \CC_{\text{id}}^n\\ & =
i_{(23)}(\CC_{(23)}^{n-1}) \bigsqcup i_{(14)}(\CC_{(14)}^{n-1})
\end{align*}
and the lassos fall into two cases:
\begin{description}
\item[Case a] Lassos with heads at
$i_{(23)}(\MB{C}_{(23)}^{n-1})$.
\end{description}
 In this case, by the induction hypothesis, one needs only to check 
lassos with head
labeled by $\Gb_{n-2}$. After writing 
 $$\CC_{(23)}^{n-1}=i_{(23)}(\CC_{\text{id}}^{n-2}
\sqcup\{\tilde\rho_\emptyset^{n-2}\})\bigsqcup
i_{(14)}(\CC_{(14)(23)}^{n-2})$$
one observes that those lassos actually reside in 
$i_{(23)}i_{(23)}(\MB{C}_{\text{id}}^{n-2}\sqcup\{\tilde\rho_\emptyset^{n-2}\})$
and those lassos with heads at a vertex of 
$i_{(23)}i_{(23)}(\MB{C}_{\text{id}}^{n-2})$ labeled by $\Gb_{n-2}$,
 are equal to $\Gb_{n-2}$. So one needs only to check the lasso 
with head at   $\tilde\rho_{n-2,n-1}^n$.\n
If $n=8$ use the path $w= \Gb_5\delta_4\Gb_6\Gb_5$ to connect 
 $\rho_{56}^8$ to $\tilde\rho_{56}^8$ and observe that 
$w\Gb_6w^{-1}=\delta_4$. Thus the remaining lasso is equal 
to  to the lasso with head at  $\rho_{56}^8$ labeled by $\Gd_4$.
This is by definition equal to $\Gd_6$ (since the path 
$\Gb_3\Gb_4\Gb_5\Gb_6\Gb_2\Gb_3\Gb_4\Gb_5$ connects $\rho_{23}^8$
 to $\rho_{56}^8$ ).

For the  $n\geq 10$ case we need to consider the following braid 
which is a generalization of $\Gd_4$ and $\Gd_6$: for even $n\geq 4$ let 
$$\Gd_n:=
[\Gb_n]\Gb_{n-1}\Gb_{n-2}\Gb_{n-3}\cdots\Gb_1^2\Gb_2^{-1}\cdots\Gb_{n-4}^{-1}
\Gb_{n-3}\Gb_{n-2}^2\Gb_{n-3}\cdots\Gb_1 .$$
\begin{lem}
  $\Gd_{n-2}$ belongs to $G$.
\end{lem}
\begin{proof}
   Notice that for $n-2\geq 6$, $\delta_{n-2}$ fixes
pointwise a small disc containing $A_2$ and $A_3$ but no other
 branch value. Connect this disc to $\partial D^2$ by an arc that
 misses the points $A_i$, the loops $\alpha_i$ and the intervals
$x_i$. Remove the disc and the arc to get a {\em 3--fold}
covering of a disc which is equivalent to $\rho(n-2)$ via an
equivalence that sends $d_6$ to $d_4$. Therefore $\Gd_{n-2}$ is a word 
in the BW generators. Glue the disc and the arc back in to get an
expression of $\delta_{n-2}$ as a word in the generators of $G$.
\end{proof}
Returning to the  $n\geq 10$ case,
 the path 
$w=\delta_{n-2}\Gb_{n-4}\Gb_{n-5}\Gb_{n-3}\Gb_{n-4}\Gb_{n-2}\Gb_{n-3}$
connects   
$\tilde\rho_{n-2,n-1} $ to $\tilde\rho_{n-5,n-3,n-2,n-1}^n 
\in \tilde\rho_{n-5,n-3,n-2,n-1}^n$.
Now observe that 
$ w\Gb_{n-2}w^{-1}=\delta_{n-2}\Gb_{n-5}\delta_{n-2}^{-1}.$ 
Therefore, by the induction hypothesis, the lasso
with head at $\tilde\rho_{n-2,n-1}$ labeled by
$\Gb_{n-2}$ is indeed in $G$.
\begin{description}
\item[Case b] Lassos with heads at $i_{(14)}(\MB{C}_{(14)}^{n-1})$.
\end{description}
 By symmetry and  Case~a) one needs to check only
lassos with head at $\tilde\rho_{23}^n$ and when $n\geq10$ the
lasso with head at $\tilde\rho_{56}^n$ labeled by
$\Gb_{n-2}$ (corresponding to $\delta_6$).\n
For the lassos based at $\tilde\rho_{23}^n$ and head labeled by
$\Gb_2$, $\Gb_4$, $\Gb_5,\ldots$, $\Gb_{n-2}$ notice that the path
$w=\Gb_{n-3}\Gb_{n-4}\dotsm\Gb_2\Gb_{n-2}\Gb_{n-3}\dotsm\Gb_3$ connects
$\tilde\rho_{n-2,n-1}^n$ to $\tilde\rho_{23}^n$. Straightforward but 
tedious calculations then give: 
 $$ w\Gb_2w^{-1}=\Gb_{n-2}\quad \text{and,}$$
  $$w\Gb_iw^{-1}=\Gb_{i-2}\quad \text{for} \quad i=4,\ldots,n-2.$$
All these lassos therefore belong to $G$.\n
 For the lasso with head at $\tilde\rho_{56}^n$ labeled by
$\Gb_{n-2}$ notice that the path $w=\Gb_{n-2}\Gb_{n-3}\dotsm\Gb_7$
connects $\tilde\rho_{5,n-1}^n$ to $\tilde\rho_{56}^n$ and that
$w\Gb_{n-2}w^{-1}=\Gb_{n-3}$ and  thus, since $\tilde\rho_{5,n-1}^n\in
i_{(23)}(\MB{C}_{(23)}^{n-1})$  this lasso is in $G$ by Case~a).\n
This concludes the proof for $\overline{\Au}(\rho_{23}^n)$ when $n$
is even.

Continuing with $n$ even:
 \begin{align*} \CC(\rho_{2,3,n-1}^n)&=\CC_{(23)(14)}^n\\
                  &=i_{(23)}(\CC_{(14)}^{n-1}\sqcup \{\rho_\emptyset^{n-1}\})
                  \bigsqcup i_{(14)}(\CC_{(23)}^{n-1}\sqcup 
                    \{\tilde\rho_\emptyset^{n-1}\})
\end{align*}
and the lassos fall into four cases:
\begin{description}
\item[Case a] Lassos with heads at
$i_{(14)}(\MB{C}_{(23)}^{n-1})$.
\end{description}
In this case, by the induction hypothesis, one only needs to check
lassos with
 head labeled by $\Gb_{n-2}$. But these heads actually lie inside
 $i_{(14)}i_{(14)}(\MB{C}_{(14)(23)}^{n-2})$, which is connected and all of its edges
commute with $\Gb_{n-2}$. It therefore follows that all of these
lassos
 are equal to the lasso
 with head at $\rho_{2,n-2,n-1}$ labeled by $\Gb_{n-2}$. Noticing that the path
 $w=\Gb_3\Gb_4\dotsm\Gb_{n-3}$ connects $\rho_{2,3,n-1}^n$ to $\rho_{2,n-2,n-1}^n$,
 one sees that the later lasso equals $[\Gb_{n-2}]\Gb_{n-3}^{-1}\dotsm\Gb_3^{-1}$.
\begin{description}
\item[Case b] Lassos with heads at $\rho_{n-1}^n$.
\end{description}
Using as tail the path
$w=\Gb_{n-2}\Gb_{n-3}\dotsm\Gb_5\delta_4\Gb_4\Gb_5\dotsm\Gb_{n-2}$, one
(straightforwardly, but tediously) calculates:
$$w\Gb_iw^{-1}= \Gb_i \quad  \text{for}\quad i\neq 3,4\, ,$$
$$w\Gb_4w^{-1}=\Gd_4 \,,$$
$$w\Gd_4w^{-1}=\Gb_4 \,.$$
Finally,
\begin{align*}w\Gb_3w^{-1}&=[\Gb_3]\Gb_{n-2}^{-1}\Gb_{n-3}^{-1}\dotsm\Gb_5^{-1}
         \Gb_4^{-1}\delta_4^{-1}\Gb_5^{-1}\dotsm\Gb_{n-3}^{-1}\Gb_{n-2}^{-1}\\
      &=[\Gb_3]\Gb_4^{-1}\delta_4^{-1}\Gb_5^{-1}\dotsm\Gb_{n-3}^{-1}\Gb_{n-2}^{-1}\\
      &=[\Gb_3]\delta_4^{-1}\Gb_4^{-1}\Gb_5^{-1}\dotsm\Gb_{n-3}^{-1}\Gb_{n-2}^{-1}\\ 
       &=[\Gb_3]\Gb_4^{-1}\Gb_5^{-1}\dotsm\Gb_{n-3}^{-1}\Gb_{n-2}^{-1}
 \end{align*}
which, as seen in Figure \ref{fig:handle}, equals
$[\Gb_{n-2}]\Gb_{n-3}^{-1}\Gb_{n-4}^{-1}\dotsm\Gb_4^{-1}$.
\begin{figure}[htp]
\centering
\mbox{\scalebox{.55}{
\begin{pspicture}(1,3.2)(7,9.6)
 \psdots(1,9)(2,9)(3,9)(5,9)(6,9)(7,9)
        (1,4)(2,4)(3,4)(5,4)(6,4)(7,4)
        \psdots[dotsize=0.05](3.4,9)(3.8,9)(4.2,9)(4.6,9)
                            (3.4,4)(3.8,4)(4.2,4)(4.6,4)
 \rput(1,9.3){14}
 \rput(2,9.3){23}
 \rput(3,9.3){23}
                 \rput(4,9.3){23's}
 \rput(5,9.3){23}
 \rput(6,9.3){23}
 \rput(7,9.3){14}
 \rput(1,3.7){14}
 \rput(2,3.7){23}
 \rput(3,3.7){23}
                  \rput(4,3.7){23's}
 \rput(5,3.7){23}
 \rput(6,3.7){23}
 \rput(7,3.7){14}
 \psline(1,7)(1.9,6.1)
 \psline(1.9,6.9)(1.56,6.56)
 \psline(1.44,6.44)(1,6)
 \psframe(2.7,4.6)(5.3,8.4)
 \psline(1,6)(1,4)
 \psline(1,9)(1,7)
 \psline(3,9)(3,8.4)
 \psline(3,4.6)(3,4)
 \psline(5,9)(5,8.4)
 \psline(5,4.6)(5,4)
 \psline(2,9)(2,4)
 \psline(6,9)(6,4)
  \psline(2.1,7.04)(2.6,7.24)
   \psline(2.8,7.32)(5.2,8.28)
   \psline(5.4,8.36)(5.9,8.56)
    \psline(6.1,8.64)(7,9)
    \psline(2.1,5.96)(2.6,5.76)
    \psline(2.8,5.68)(5.2,4.72)
    \psline(5.4,4.64)(5.9,4.44)
    \psline(6.1,4.36)(7,4)
                          \psline[linewidth=.005](3.8,9)(3.8,8.4)
                          \psline[linewidth=.005](4.2,9)(4.2,8.4)
                          \psline[linewidth=.005](4.6,9)(4.6,8.4)
                          \psline[linewidth=.005](3.4,9)(3.4,8.4)
                        \psline[linewidth=.005](3.8,4.6)(3.8,4)
                        \psline[linewidth=.005](4.2,4.6)(4.2,4)
                        \psline[linewidth=.005](4.6,4.6)(4.6,4)
                        \psline[linewidth=.005](3.4,4.6)(3.4,4)
                 \rput(3.5,5){id}
\end{pspicture}
\hspace{2cm}
\begin{pspicture}(1,3.2)(7,9.6)
 \psdots(1,9)(2,9)(3,9)(5,9)(6,9)(7,9)
        (1,4)(2,4)(3,4)(5,4)(6,4)(7,4)
        \psdots[dotsize=0.05](3.4,9)(3.8,9)(4.2,9)(4.6,9)
                            (3.4,4)(3.8,4)(4.2,4)(4.6,4)
 \rput(1,9.3){14}
 \rput(2,9.3){23}
 \rput(3,9.3){23}
                 \rput(4,9.3){23's}
 \rput(5,9.3){23}
 \rput(6,9.3){23}
 \rput(7,9.3){14}
 \rput(1,3.7){14}
 \rput(2,3.7){23}
 \rput(3,3.7){23}
                  \rput(4,3.7){23's}
 \rput(5,3.7){23}
 \rput(6,3.7){23}
 \rput(7,3.7){14}
  \psline(6.1,6.9)(7,6)
  \psline(7,7)(6.56,6.56)
  \psline(6.44,6.44)(6.1,6.1)
 \psframe(2.7,4.6)(5.3,8.4)
 \psline(7,6)(7,4)
 \psline(7,9)(7,7)
 \psline(3,9)(3,8.4)
 \psline(3,4.6)(3,4)
 \psline(5,9)(5,8.4)
 \psline(5,4.6)(5,4)
 \psline(2,9)(2,4)
 \psline(6,9)(6,4)
    \psline(5.9,5.96)(5.4,5.76)
    \psline(5.2,5.68)(2.8,4.72)
    \psline(2.6,4.64)(2.1,4.44)
    \psline(1.9,4.36)(1,4)
       \psline(5.9,7.04)(5.4,7.24)
       \psline(5.2,7.32)(2.8,8.28)
       \psline(2.6,8.36)(2.1,8.56)
       \psline(1.9,8.64)(1,9)
                         \psline[linewidth=.005](3.8,9)(3.8,8.4)
                          \psline[linewidth=.005](4.2,9)(4.2,8.4)
                          \psline[linewidth=.005](4.6,9)(4.6,8.4)
                          \psline[linewidth=.005](3.4,9)(3.4,8.4)
                        \psline[linewidth=.005](3.8,4.6)(3.8,4)
                        \psline[linewidth=.005](4.2,4.6)(4.2,4)
                        \psline[linewidth=.005](4.6,4.6)(4.6,4)
                        \psline[linewidth=.005](3.4,4.6)(3.4,4)
                 \rput(5,5){id}
\end{pspicture}
\hspace{2cm}
\begin{pspicture}(1,3.2)(7,9.6)
 \psdots(1,9)(2,9)(3,9)(5,9)(6,9)(7,9)
        (1,4)(2,4)(3,4)(5,4)(6,4)(7,4)
        \psdots[dotsize=0.05](3.4,9)(3.8,9)(4.2,9)(4.6,9)
                            (3.4,4)(3.8,4)(4.2,4)(4.6,4)
 \rput(1,9.3){14}
 \rput(2,9.3){23}
 \rput(3,9.3){23}
                 \rput(4,9.3){23's}
 \rput(5,9.3){23}
 \rput(6,9.3){23}
 \rput(7,9.3){14}
 \rput(1,3.7){14}
 \rput(2,3.7){23}
 \rput(3,3.7){23}
                  \rput(4,3.7){23's}
 \rput(5,3.7){23}
 \rput(6,3.7){23}
 \rput(7,3.7){14}
  \psline(6,7)(7,6)
  \psline(7,7)(6.56,6.56)
  \psline(6.44,6.44)(6,6)
 \psline(2.7,4.6)(5.3,4.6)
 \psline(2.7,8.4)(5.3,8.4)
 \psline(2.7,4.6)(2.7,4.63)
 \psline(2.7,4.8)(2.7,8.25)
 \psline(5.3,4.6)(5.3,5.6)
 \psline(5.3,5.8)(5.3,7.2)
 \psline(5.3,7.4)(5.3,8.4)
 \psline(7,6)(7,4)
 \psline(7,9)(7,7)
 \psline(3,9)(3,8.4)
 \psline(3,4.6)(3,4)
 \psline(5,9)(5,8.4)
 \psline(5,4.6)(5,4)
  \psline(2,9)(2,8.7)
  \psline(2,8.5)(2,4.5)
  \psline(2,4.3)(2,4)
  \psline(6,9)(6,7.1)
  \psline(6,6.9)(6,6.1)
  \psline(6,5.9)(6,4)
    \psline(6,6)(1,4)
       \psline(6,7)(1,9)
                         \psline[linewidth=.005](3.8,9)(3.8,8.4)
                          \psline[linewidth=.005](4.2,9)(4.2,8.4)
                          \psline[linewidth=.005](4.6,9)(4.6,8.4)
                          \psline[linewidth=.005](3.4,9)(3.4,8.4)
                        \psline[linewidth=.005](3.8,4.6)(3.8,4)
                        \psline[linewidth=.005](4.2,4.6)(4.2,4)
                        \psline[linewidth=.005](4.6,4.6)(4.6,4)
                        \psline[linewidth=.005](3.4,4.6)(3.4,4)
                 \rput(5,5){id}
\end{pspicture}
}}
 \caption{First shift to the right by an isotopy
and then use move $\mathcal{P}$ to pass over the strands in the
box.}\label{fig:handle}
\end{figure}
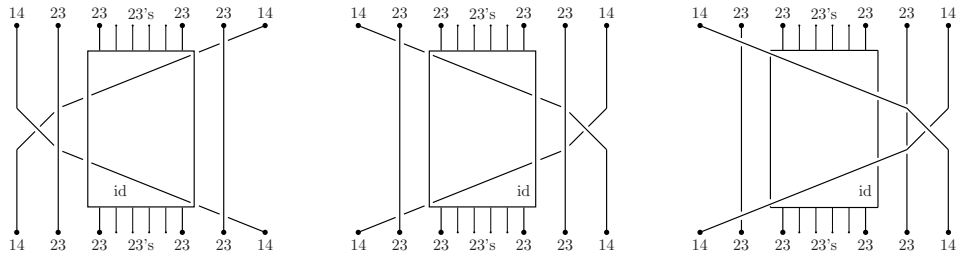

\begin{description}
\item[Case c] Lassoes with heads at
$i_{(23)}(\MB{C}_{(14)}^{n-1})$.
\end{description}
By symmetry and Case~a) one needs only to check lassoes with head
at $\tilde\rho_{2,3,n-1}^n$, the lasso with head at
$\tilde\rho_{2,n-2,n-1}$ labelled by $\Gb_{n-2}$ (which corresponds to 
$[\Gb_{n-2}]\Gb_{n-3}^{-1}\dotsm\Gb_3^{-1}$), and for $n\geq 10$ the
lasso with head at $\tilde\rho_{2,3,n-1}^n$ labelled by
$\delta_4$ (which corresponds to $\delta_6$).\n
For the lassos with head at $\tilde\rho_{2,3,n-1}^n$ labelled by
$\delta_4$ or $\Gb_i$ for $i=2,4,5,\ldots,n-4$ notice that
$\Gb_{n-2}$ connects $\tilde\rho_{2,3,n-1}^n$ to
$\tilde\rho_{2,3,n-2}^n\in i_{(14)}(\MB{C}_{(12)}^{n-1})$ and so they
have been checked in Case~a). For the one with head labelled by
$\Gb_{n-3}$ ``pull further into $i_{(14)}(\MB{C}_{(12)}^{n-1})$''
using the path $\Gb_{n-2}\Gb_{n-3}$.\n
For the lasso with head at $\tilde\rho_{2,n-2,n-1}$ labelled by
$\Gb_{n-2}$ use again the path $\Gb_{n-2}\Gb_{n-3}$ to ``pull it'' to
the lasso with head at $\tilde\rho_{2,n-3,n-2}\in
i_{(14)}(\MB{C}_{(12)}^{n-1})$, labelled
by $\Gb_{n-3}$.\n
Finally for the lasso with head at $\tilde\rho_{2,3,n-1}^n$
labelled by $\delta_4$ use $\Gb_{n-2}$ to ``pull'' it to
$\tilde\rho_{2,3,n-1}^n \in i_{(14)}(\MB{C}_{(23)}^{n-1})$ with label
$\delta_4$.
\begin{description}
\item[Case d] Lassoes with heads at $\tilde\rho_{n-1}^n$.
\end{description}
Analogously to Case~b), the check for these lassos 
reduces to Case~c).\n
This concludes the proof for $\overline{\Au}(\rho_{2,3,n-1}^n)$ and
 the even $n$ case.
\end{proof}

\section{The kernel}
\label{sec:kernel}
\begin{figure}[ht!]
\centerline{\noindent\hglue 2cm
\mbox{\scalebox{.83}{
\begin{pspicture}(-1,-1)(15,15)
 \psellipse[linecolor=blue](.5,1.5)(.5,.3)
 \psellipse[linecolor=blue](2.5,1.5)(.5,.3)
 \psellipse[linecolor=blue](4.5,1.5)(.5,.3)
 \psellipse[linecolor=blue](6.5,1.5)(.5,.3)
 \psellipse[linecolor=blue](8.5,1.5)(.5,.3)
 \psellipse[linecolor=blue](12.5,1.5)(.5,.3)
 \psbezier(1,1.5)(1.25,1)(1.75,1)(2,1.5)
  \psbezier(3,1.5)(3.25,1)(3.75,1)(4,1.5)
   \psbezier(5,1.5)(5.25,1)(5.75,1)(6,1.5)
    \psbezier(7,1.5)(7.25,1)(7.75,1)(8,1.5)
     \psbezier(9,1.5)(9.25,1)(9.75,1)(10,1.5)
          \psbezier(11,1.5)(11.25,1)(11.75,1)(12,1.5)
     \psbezier
     (0,1.5)(-.5,-1)(13.5,-1)(13,1.5)
  \psframe[linecolor=white,fillstyle=solid,fillcolor=white](9.5,-1)(11.5,2)
  \psdots[dotsize=.1](10.3333,.5)(10.6666,.5)(10.9999,.5)
  \psbezier[linestyle=dotted,linecolor=red](5,1.5)(5.2,0)(.8,-.8)(.5,1.8)
  \psbezier[linecolor=red](.4,1.2)(0,0)(4,.6)(4.5,1.2)
  \psbezier[linestyle=dotted,linecolor=red](4.5,1.8)(4,.6)(.8,.4)(1,1.5)
 \rput(0,6){\psellipse[linecolor=blue](2.5,1.5)(.5,.3)
 \psframe[linecolor=white,fillstyle=solid,fillcolor=white](0,1.5)(13,2)
\psellipse[linecolor=blue](.5,1.5)(.5,.3)
 \psellipse[linestyle=dotted,linecolor=blue](2.5,1.5)(.5,.3)
  \psellipse[linecolor=blue](4.5,1.5)(.5,.3)
 \psellipse[linecolor=blue](6.5,1.5)(.5,.3)
 \psellipse[linecolor=blue](8.5,1.5)(.5,.3)
 \psellipse[linecolor=blue](12.5,1.5)(.5,.3)
 \psbezier(1,1.5)(1.25,1)(1.75,1)(2,1.5)
  \psbezier(3,1.5)(3.25,1)(3.75,1)(4,1.5)
   \psbezier(5,1.5)(5.25,1)(5.75,1)(6,1.5)
    \psbezier(7,1.5)(7.25,1)(7.75,1)(8,1.5)
     \psbezier(9,1.5)(9.25,1)(9.75,1)(10,1.5)
          \psbezier(11,1.5)(11.25,1)(11.75,1)(12,1.5)
     \psbezier
     (0,1.5)(-.5,-1)(13.5,-1)(13,1.5)
  \psframe[linecolor=white,fillstyle=solid,fillcolor=white](9.5,-1)(11.5,2)
  \psdots[dotsize=.1](10.3333,.5)(10.6666,.5)(10.9999,.5)
  \psbezier[linestyle=dotted,linecolor=red](5,1.5)(5.2,0)(.8,-.8)(.5,1.8)
  \psbezier[linecolor=red](.4,1.2)(0,0)(4,.6)(4.5,1.2)
  \psbezier[linestyle=dotted,linecolor=red](4.5,1.8)(4,.6)(.8,.4)(1,1.5)}
  \rput(0,8){\psellipse[linecolor=blue](.5,1.5)(.5,.3)
  \psframe[linecolor=white,fillstyle=solid,fillcolor=white](0,1.5)(13,2)
  \psellipse[linestyle=dotted,linecolor=blue](.5,1.5)(.5,.3)
 \psellipse[linecolor=blue](2.5,1.5)(.5,.3)
 \psellipse[linecolor=blue](4.5,1.5)(.5,.3)
 \psellipse[linecolor=blue](6.5,1.5)(.5,.3)
 \psellipse[linecolor=blue](8.5,1.5)(.5,.3)
 \psellipse[linecolor=blue](12.5,1.5)(.5,.3)
 \psbezier(1,1.5)(1.25,1)(1.75,1)(2,1.5)
  \psbezier(3,1.5)(3.25,1)(3.75,1)(4,1.5)
   \psbezier(5,1.5)(5.25,1)(5.75,1)(6,1.5)
    \psbezier(7,1.5)(7.25,1)(7.75,1)(8,1.5)
     \psbezier(9,1.5)(9.25,1)(9.75,1)(10,1.5)
         \psbezier(11,1.5)(11.25,1)(11.75,1)(12,1.5)
       \psframe[linecolor=white,fillstyle=solid,fillcolor=white](9.5,-1)(11.5,2)
  \psdots[dotsize=.1](10.3333,.5)(10.6666,.5)(10.9999,.5)}
  \rput(0,11){\psellipse[linecolor=blue](.5,1.5)(.5,.3)
  \psellipse[linecolor=blue](2.5,1.5)(.5,.3)
 \psellipse[linecolor=blue](4.5,1.5)(.5,.3)
 \psellipse[linecolor=blue](6.5,1.5)(.5,.3)
 \psellipse[linecolor=blue](8.5,1.5)(.5,.3)
 \psellipse[linecolor=blue](12.5,1.5)(.5,.3)
 \psframe[linecolor=white,fillstyle=solid,fillcolor=white](0,1.5)(13,2)
 \psellipse[linestyle=dotted,linecolor=blue](2.5,1.5)(.5,.3)
 \psellipse[linestyle=dotted,linecolor=blue](4.5,1.5)(.5,.3)
 \psellipse[linestyle=dotted,linecolor=blue](6.5,1.5)(.5,.3)
 \psellipse[linestyle=dotted,linecolor=blue](8.5,1.5)(.5,.3)
 \psellipse[linestyle=dotted,linecolor=blue](12.5,1.5)(.5,.3)
 \psbezier(1,1.5)(1.25,1)(1.75,1)(2,1.5)
   \psbezier(5,1.5)(5.25,1)(5.75,1)(6,1.5)
    \psbezier(7,1.5)(7.25,1)(7.75,1)(8,1.5)
     \psbezier(9,1.5)(9.25,1)(9.75,1)(10,1.5)
          \psbezier(11,1.5)(11.25,1)(11.75,1)(12,1.5)
  \psframe[linecolor=white,fillstyle=solid,fillcolor=white](9.5,-1)(11.5,2)
  \psdots[dotsize=.1](10.3333,.5)(10.6666,.5)(10.9999,.5)}
 \rput(-1,7){Sheet 4}
 \rput(-1,8.5){Sheet 1}
 \rput(-1,11){Sheet 2}
 \rput(-1,13){Sheet 3}
 \psframe[linecolor=white,fillstyle=solid,fillcolor=white](1,8)(4,11)
\psframe[linecolor=white,fillstyle=solid,fillcolor=white](0,10)(2,13)
\psellipse[linecolor=blue](.5,12.5)(.5,.3)
 \psbezier(1,12.5)(1.25,13)(1.75,13)(2,12.5)
 \psbezier(5,12.5)(5.25,13)(5.75,13)(6,12.5)
 \psbezier(3,12.5)(3.25,13)(3.75,13)(4,12.5)
 \psbezier(7,12.5)(7.25,13)(7.75,13)(8,12.5)
 \psbezier(9,12.5)(9.25,13)(9.75,13)(10,12.5)
 \psbezier(11,12.5)(11.25,13)(11.75,13)(12,12.5)
 \psbezier(0,12.5)(-.5,15)(13.5,15)(13,12.5)
 \psframe[linecolor=white,fillstyle=solid,fillcolor=white](9.5,12.5)(11.5,15)
  \psdots[dotsize=.1](10.3333,13.5)(10.6666,13.5)(10.9999,13.5)
\psbezier
     (0,9.5)(-.2,8.4)(1.8,8.6)(2,7.5)
  \psbezier
     (3,7.5)(3.2,9.8)(14,6.5)(13,9.5)
  \psbezier
     (1,9.5)(1.5,9)(3.5,9)(4,9.5)
  \psbezier
     (0,9.5)(-.1,12)(1.5,11.5)(2,11.8) 
  \psbezier
     (1,9.5)(.8,12)(13.5,9)(13,12.5)
  \psbezier[linestyle=dotted,linecolor=red](5,9.5)(5.2,7.8)(.2,7.6)(.6,9.8)
  \psbezier[linecolor=red](.4,9.2)(1,8.2)(4,8.2)(4.5,9.2)
  \psbezier[linestyle=dotted,linecolor=red](4.5,9.8)(4,8.6)(.8,8.4)(1,9.5)
  \psbezier[linecolor=red](5,12.5)(5.2,14.2)(.2,13.5)(.5,12.8)
  \psbezier[linecolor=red](4.7,12.2)(4.7,13.4)(.8,13.4)(1,12.5)
  \psbezier[linestyle=dotted,linecolor=red](4.5,12.8)(5.2,14.6)(.2,13.6)(.5,12.2)
  \psbezier[linestyle=dotted,linecolor=red](5,12.5)(5.2,10.6)(0,10.6)(.2,11)
  \psbezier[linecolor=red](.2,11)(.6,10)(.4,9.8)(.4,9.2)
  \psbezier[linestyle=dotted,linecolor=red](4.5,12.8)(4.8,11)(0,11.2)(.6,11.4)
  \psbezier[linecolor=red](.6,11.4)(.8,10)(1,9.8)(1,9.5)
  \psbezier[linestyle=dotted,linecolor=red](.6,9.8)(.8,10.4)(1.2,10.3)(1.4,10.2)
  \psbezier[linecolor=red](1.4,10.2)(.5,11.5)(5.5,10)(4.7,12.2)
\psframe[linecolor=white,fillstyle=solid,fillcolor=white](9.5,8)(11.5,8.3)
\psframe[linecolor=white,fillstyle=solid,fillcolor=white](9.5,10)(11.5,11.3)
 \psline{->}(6.5,5)(6.5,2.5)
\psellipse[linecolor=blue](2.5,11.8)(.5,.3)
  \pscurve(3,11.8)(3.6,12)(4,12.5)
\end{pspicture}
}}}
\caption{The covering $\rho$ and how $d_4$ lifts}\label{fig:model}
\end{figure}
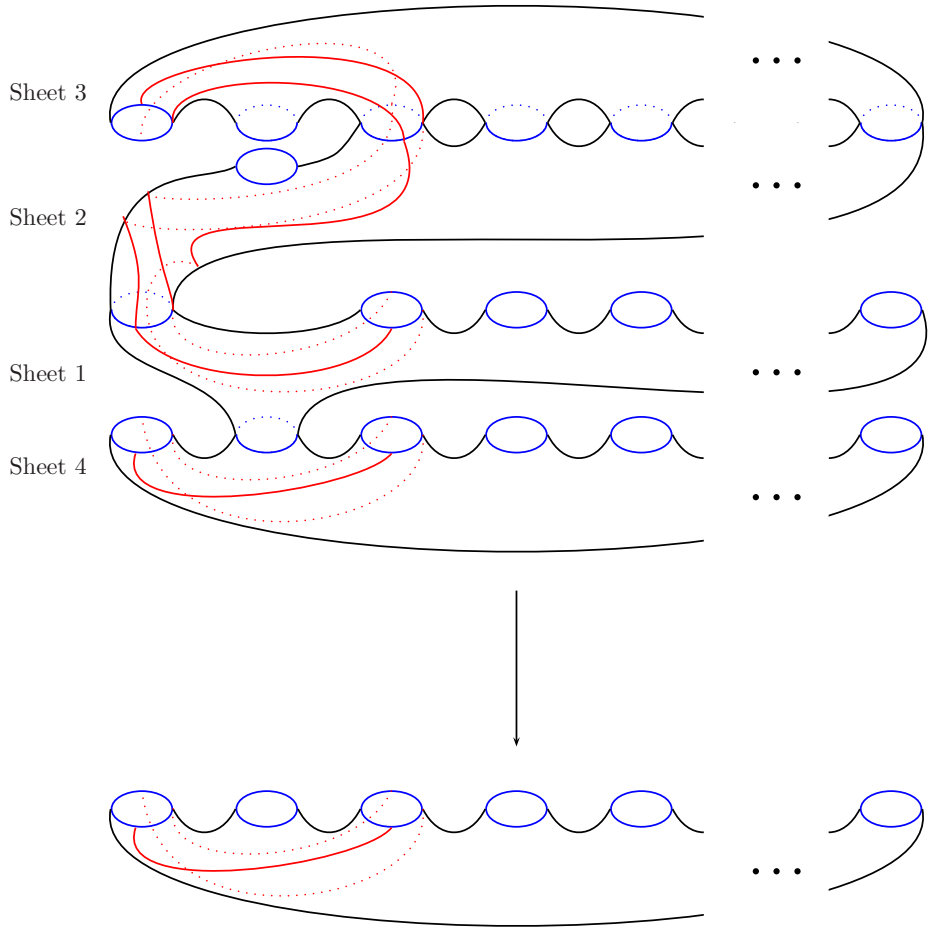

In this section we compute the kernel of the (reduced) lifting homomorphism
for the standard $4$--fold covering $\Gr^n_{23}$, which, for simplicity, 
will be denoted by $\Gr$. The genus of the total space of $\Gr$  
according to the  Riemann-Hurwitz formula (see~\cite{BernEdm1979})
is given by $g=\frac{n-6}{2}.$ 
\begin{thm}\label{thm:map}
The reduced lifting homomorphism 
$$\bar{\Gl}\co \overline{\Au}(\Gr)\to \FM_g$$
is given on the generators of Theorem~\ref{thm:gen} by  
\begin{equation}\label{eq:themap}
\bar\lambda(x)=
\begin{cases} \text{id} & \text{if} \quad x=\Gb_0,\Gb_2\\
                       a_{i-1} & \text{if $ \quad x=\Gb_{2i}$ for $i\geq2$}\\
                       b_{i-1} &\text{if  $\quad x=\Gb_{2i+1}$ with
                       $i\geq2$}\\
                       a_1 & \text{if $ \quad x=\delta_4$}\\
                       d & \text{if $ \quad x=\delta_6$}
                   \end{cases}
\end{equation}
where $a_i$, $b_i$, and $d$ are Wajnryb's generators of $\CM_g$, (see~\cite{Wajn1983}).
\end{thm}
\begin{proof}
  Since all the generators given in Theorem~\ref{thm:gen} are rotations 
around intervals, their
image under $\bar\Gl$ can be determined by lifting those intervals 
(see~\cite{BernEdm1979}).
Figure~\ref{fig:model} depicts a model of $\Gr$ constructed by 
``cutting and pasting,''
 the 
$2$--sphere $S^2$ is cut open along intervals
 $x_i$ for even $i$, and then four copies of the 
cut sphere  are glued together according to the specified monodromies. 
 As an example, in Figure~\ref{fig:model} the interval 
$\Gd_4$ is shown  and how it lifts to the disjoint union
 of a curve isotopic to $a_1$ and two arcs.
\end{proof}
Now one can prove the following
\begin{thm}\label{thm:kernel}
 The kernel of $\bar{\lambda}$ is the
smallest normal subgroup of $\overline{\Au}(\rho_{2g+6})$
containing the elements $\Gb_0$, $\Gb_2$, $B$, $D$, and all words
obtained from $B$ or $D$ by replacing some appearances of $\Gb_4$
with $\delta_4$, where:
$$B=(\Gb_4\Gb_5\Gb_6)^4([\delta_6^{-1}]\Gb_6^{-1}\Gb_5^{-1}\Gb_4^{-2}
\Gb_5^{-1}\Gb_6^{-1}\Gb_7^{-1})\delta_6^{-1},$$
$$D=\Gb_{2g+4}\chi\Gb_{2g+4}^{-1}\chi^{-1},$$
where:
$$\chi=\Gb_{2g+3}\Gb_{2g+2}\dotsm\Gb_5\Gb_4^2\Gb_5\dotsm\Gb_{2g+2}\Gb_{2g+3},$$
\end{thm}
\begin{proof} Given the Wajnryb's presentation of $\CM_g$  and
the fact that $\bar{\lambda}$ is given by \eqref{eq:themap} the proof
reduces to checking that the words obtained by the relators in~\cite{Wajn1983}, 
after replacing  the generators of
$\CM_g$ with their preimages to generators of
$\overline{\Au}(\rho_{2g+6})$ in all possible ways, are in the normal closure of the
given elements.\n
 The calculations are identical to those in Theorems 5.1 and 6.1 of
 \cite{BirWajn1985} after replacing  $x_i$ in \cite{BirWajn1985} with $\Gb_{i+2}$
 for $i\geq2$, and $\delta_4$ in \cite{BirWajn1985} with $\delta_6$,
 and then replacing $\Gb_4$ with $\delta_4$ in all possible ways.
\end{proof}

\section{Moves}
\label{sec:moves}

In this section we prove the  two theorems mentioned in the introduction.
Recall that ``bi-tricolored'' means ``colored by transpositions of $\CS_4$.''

\begin{figure}[ht!]\vspace{-.3cm}
\centerline{\noindent\hglue 1.2cm\mbox{
      \subfigure{\scalebox{0.8}{
\begin{pspicture}(.5,0)(4,5)
 \psdots(0,4)(.5,4)(1,4)(1.5,4)(2,4)(2.5,4)(4,4)(4.5,4)
        (0,1)(.5,1)(1,1)(1.5,1)(2,1)(2.5,1)(4,1)(4.5,1)
   \psline(0,4)(0,1)
   \psline(0.5,4)(0.5,1)
   \psline(2,4)(2,1)
   \psline(2.5,4)(2.5,1)
   \psline(4,4)(4,1)
   \psline(4.5,4)(4.5,1)
 \psline[linearc=.25](1.5,4)(1.5,2.75)(1,2.25)(1,1)
\psframe[linecolor=white,fillstyle=solid,fillcolor=white](1.2,2.4)(1.38,2.55)
\psline[linearc=.25](1,4)(1,2.75)(1.5,2.25)(1.5,1)
           \psframe(-.3,0)(4.8,1)
           \psframe(-.3,4)(4.8,5)
    \rput(0,4.3){$12$}
 \rput(0.5,4.3){$12$}
 \rput(1,4.3){$14$}
 \rput(1.5,4.3){$14$}
 \rput(2,4.3){$23$}
 \rput(2.5,4.3){$23$}
 \rput(4,4.3){$23$}
 \rput(4.5,4.3){$23$}
                   \rput(3.2,4.3){$23$'s}
 \rput(0,.7){$12$}
 \rput(0.5,.7){$12$}
 \rput(1,.7){$14$}
 \rput(1.5,.7){$14$}
 \rput(2,.7){$23$}
 \rput(2.5,.7){$23$}
 \rput(4,.7){$23$}
 \rput(4.5,.7){$23$}
                   \rput(3.2,.7){$23$'s}
  \psdots[dotsize=0.1](3,4)(3.3,4)(3.6,4)(3,1)(3.3,1)(3.6,1)
                      (3,2.5)(3.3,2.5)(3.6,2.5)
     \rput(2.5,4.8){$\tau$}
     \rput(2.5,.2){$\tau'$}
\end{pspicture} 
}
}
\quad
\begin{pspicture}(4,5)
\psline{<->}(1,2.5)(3,2.5)
\rput(2,2.8){\large{I}}
\end{pspicture}
\quad
\subfigure{\scalebox{0.8}{
\begin{pspicture}(.5,0)(6,5)
  \psdots(0,4)(.5,4)(1,4)(1.5,4)(2,4)(2.5,4)(4,4)(4.5,4)
        (0,1)(.5,1)(1,1)(1.5,1)(2,1)(2.5,1)(4,1)(4.5,1)
   \psline(0,4)(0,1)
   \psline(0.5,4)(0.5,1)
   \psline(2,4)(2,1)
   \psline(2.5,4)(2.5,1)
   \psline(4,4)(4,1)
   \psline(4.5,4)(4.5,1)
 \psline(1.5,4)(1.5,1)
\psline(1,4)(1,1)
           \psframe(-.3,0)(4.8,1)
           \psframe(-.3,4)(4.8,5)
    \rput(0,4.3){$12$}
 \rput(0.5,4.3){$12$}
 \rput(1,4.3){$14$}
 \rput(1.5,4.3){$14$}
 \rput(2,4.3){$23$}
 \rput(2.5,4.3){$23$}
 \rput(4,4.3){$23$}
 \rput(4.5,4.3){$23$}
                   \rput(3.2,4.3){$23$'s}
 \rput(0,.7){$12$}
 \rput(0.5,.7){$12$}
 \rput(1,.7){$14$}
 \rput(1.5,.7){$14$}
 \rput(2,.7){$23$}
 \rput(2.5,.7){$23$}
 \rput(4,.7){$23$}
 \rput(4.5,.7){$23$}
                   \rput(3.2,.7){$23$'s}
  \psdots[dotsize=0.1](3,4)(3.3,4)(3.6,4)(3,1)(3.3,1)(3.6,1)
                      (3,2.5)(3.3,2.5)(3.6,2.5)
     \rput(2.5,4.8){$\tau$}
     \rput(2.5,.2){$\tau'$}
\end{pspicture}
}
}
}}
\vspace{-.3cm}
 \centerline{     \mbox{
      \subfigure{\scalebox{0.8}{
\begin{pspicture}(.5,0)(4,5)
 \psdots(0,4)(.5,4)(1,4)(1.5,4)(2,4)(2.5,4)(4,4)(4.5,4)
        (0,1)(.5,1)(1,1)(1.5,1)(2,1)(2.5,1)(4,1)(4.5,1)
   \psline(1,4)(1,1)
   \psline(1.5,4)(1.5,1)
   \psline(2,4)(2,1)
   \psline(2.5,4)(2.5,1)
   \psline(4,4)(4,1)
   \psline(4.5,4)(4.5,1) 
 \psline[linearc=.25](.5,4)(.5,3)(0,2)(0,1)
\psframe[linecolor=white,fillstyle=solid,fillcolor=white](.2,2.4)(.39,2.6)
\psline[linearc=.25](0,4)(0,3)(.5,2)(.5,1)
           \psframe(-.3,1)(4.8,0)
           \psframe(-.3,4)(4.8,5)
    \rput(0,4.3){$12$}
 \rput(0.5,4.3){$12$}
 \rput(1,4.3){$14$}
 \rput(1.5,4.3){$14$}
 \rput(2,4.3){$23$}
 \rput(2.5,4.3){$23$}
 \rput(4,4.3){$23$}
 \rput(4.5,4.3){$23$}
                   \rput(3.2,4.3){$23$'s}
 \rput(0,.7){$12$}
 \rput(0.5,.7){$12$}
 \rput(1,.7){$14$}
 \rput(1.5,.7){$14$}
 \rput(2,.7){$23$}
 \rput(2.5,.7){$23$}
 \rput(4,.7){$23$}
 \rput(4.5,.7){$23$}
                   \rput(3.2,.7){$23$'s}
  \psdots[dotsize=0.1](3,4)(3.3,4)(3.6,4)(3,1)(3.3,1)(3.6,1)
                      (3,2.5)(3.3,2.5)(3.6,2.5)
     \rput(2.5,4.8){$\tau$}
     \rput(2.5,.2){$\tau'$} 
\end{pspicture}
}
}
\quad
\begin{pspicture}(4,5)
\psline{<->}(1,2.5)(3,2.5)
\rput(2,2.8){\large{II}}
\end{pspicture}
\quad
\subfigure{\scalebox{0.8}{
\begin{pspicture}(4,5)
  \psdots(0,4)(.5,4)(1,4)(1.5,4)(2,4)(2.5,4)(4,4)(4.5,4)
        (0,1)(.5,1)(1,1)(1.5,1)(2,1)(2.5,1)(4,1)(4.5,1)
   \psline(0,4)(0,1)
   \psline(0.5,4)(0.5,1)
   \psline(2,4)(2,1)
   \psline(2.5,4)(2.5,1)
   \psline(4,4)(4,1)
   \psline(4.5,4)(4.5,1)
 \psline(1.5,4)(1.5,1)
\psline(1,4)(1,1)
           \psframe(-.3,0)(4.8,1)
           \psframe(-.3,4)(4.8,5)
    \rput(0,4.3){$12$}
 \rput(0.5,4.3){$12$}
 \rput(1,4.3){$14$}
 \rput(1.5,4.3){$14$}
 \rput(2,4.3){$23$}
 \rput(2.5,4.3){$23$}
 \rput(4,4.3){$23$}
 \rput(4.5,4.3){$23$}
                   \rput(3.2,4.3){$23$'s}
 \rput(0,.7){$12$}
 \rput(0.5,.7){$12$}
 \rput(1,.7){$14$}
 \rput(1.5,.7){$14$}
 \rput(2,.7){$23$}
 \rput(2.5,.7){$23$}
 \rput(4,.7){$23$}
 \rput(4.5,.7){$23$}
                   \rput(3.2,.7){$23$'s}
  \psdots[dotsize=0.1](3,4)(3.3,4)(3.6,4)(3,1)(3.3,1)(3.6,1)
                      (3,2.5)(3.3,2.5)(3.6,2.5)
     \rput(2.5,4.8){$\tau$}
     \rput(2.5,.2){$\tau'$}
\end{pspicture}
}
}
}}
\vspace{-0.8cm}
\centerline{\noindent\hglue 1.4cm \mbox{
      \subfigure{\scalebox{0.7}{
\begin{pspicture}(1.5,0)(5,5)
 \psdots(0,4)(.5,4)(1,4)(1.5,4)(2,4)(2.5,4)(3,4)(3.5,4)(5.5,4)(6,4)
        (0,1)(.5,1)(1,1)(1.5,1)(2,1)(2.5,1)(3,1)(3.5,1)(5.5,1)(6,1)
   \psline(0,4)(0,1)
   \psline(0.5,4)(0.5,1)
   \psline(3,4)(3,1)
   \psline(3.5,4)(3.5,1)
   \psline(5.5,4)(5.5,1)
   \psline(6,4)(6,1)
  \psline[linearc=.25](1,4)(2,2.5)(1.8,2.2)
  \psline[linearc=.25](1.5,4)(2.5,2.5)(2.1,1.9)
  \psline[linearc=.25](1.4,3.1)(1,2.5)(2,1)
  \psline[linearc=.25](1.7,2.8)(1.5,2.5)(2.5,1)
   \psline(1.6,3.4)(1.7,3.55)
   \psline(1.85,3.755)(2,4)
    \psline(1.85,3.025)(1.95,3.175)
    \psline(2.05,3.325)(2.5,4)
    \psline(1.65,1.975)(1.56,1.84)
    \psline(1.4,1.6)(1,1)
    \psline(1.9,1.6)(1.8,1.45)
    \psline(1.7,1.3)(1.5,1)
           \psframe(-.3,0)(6.3,1)
           \psframe(-.3,4)(6.3,5)
 \rput(0,4.3){$14$}
 \rput(0.5,4.3){$14$}
 \rput(1,4.3){$12$}
 \rput(1.5,4.3){$12$}
 \rput(2,4.3){$23$}
 \rput(2.5,4.3){$23$}
 \rput(3,4.3){$23$}
 \rput(3.5,4.3){$23$}
 \rput(5.5,4.3){$23$}
 \rput(6,4.3){$23$}
                   \rput(4.35,4.3){$23$'s}
 \rput(0,.7){$14$}
 \rput(0.5,.7){$14$}
 \rput(1,.7){$12$}
 \rput(1.5,.7){$12$}
 \rput(2,.7){$23$}
 \rput(2.5,.7){$23$}
 \rput(3,.7){$23$}
 \rput(3.5,.7){$23$}
 \rput(5.5,.7){$23$}
 \rput(6,.7){$23$}
                   \rput(4.35,.7){$23$'s}
  \psdots[dotsize=0.1](4.2,4)(4.5,4)(4.8,4)(4.2,1)(4.5,1)(4.8,1)
                      (4.2,2.5)(4.5,2.5)(4.8,2.5)
     \rput(3.2,4.8){$\tau$}
     \rput(3.2,.2){$\tau'$}
\end{pspicture}
}
}
\quad
\begin{pspicture}(4,5)
\psline{<->}(1,2)(3,2)
\rput(2,2.3){\large{III}}
\end{pspicture}
\quad
\subfigure{\scalebox{0.7}{
\begin{pspicture}(.5,0)(6.5,5)
 \psdots(0,4)(.5,4)(1,4)(1.5,4)(2,4)(2.5,4)(3,4)(3.5,4)(5.5,4)(6,4)
        (0,1)(.5,1)(1,1)(1.5,1)(2,1)(2.5,1)(3,1)(3.5,1)(5.5,1)(6,1)
           \psframe(-.3,0)(6.3,1)
           \psframe(-.3,4)(6.3,5)
 \rput(0,4.3){$14$}
 \rput(0.5,4.3){$14$}
 \rput(1,4.3){$12$}
 \rput(1.5,4.3){$12$}
 \rput(2,4.3){$23$}
 \rput(2.5,4.3){$23$}
 \rput(3,4.3){$23$}
 \rput(3.5,4.3){$23$}
 \rput(5.5,4.3){$23$}
 \rput(6,4.3){$23$}
                   \rput(4.35,4.3){$23$'s}
 \rput(0,.7){$14$}
 \rput(0.5,.7){$14$}
 \rput(1,.7){$12$}
 \rput(1.5,.7){$12$}
 \rput(2,.7){$23$}
 \rput(2.5,.7){$23$}
 \rput(3,.7){$23$}
 \rput(3.5,.7){$23$}
 \rput(5.5,.7){$23$}
 \rput(6,.7){$23$}
                   \rput(4.35,.7){$23$'s}
  \psdots[dotsize=0.1](4.2,4)(4.5,4)(4.8,4)(4.2,1)(4.5,1)(4.8,1)
                      (4.2,2.5)(4.5,2.5)(4.8,2.5)
    \rput(3.2,4.8){$\tau$}
     \rput(3.2,.2){$\tau'$}
 \psline(0,4)(0,1)
   \psline(0.5,4)(0.5,1)
   \psline(3,4)(3,1)
   \psline(3.5,4)(3.5,1)
   \psline(5.5,4)(5.5,1)
   \psline(6,4)(6,1)
  \psline(1,4)(1,1)
  \psline(1.5,4)(1.5,1)
   \psline(2,4)(2,1)
      \psline(2.5,4)(2.5,1)
\end{pspicture}
}
}
}}
\caption{The first three non-local moves}
 \label{fig:movesI}
\end{figure}
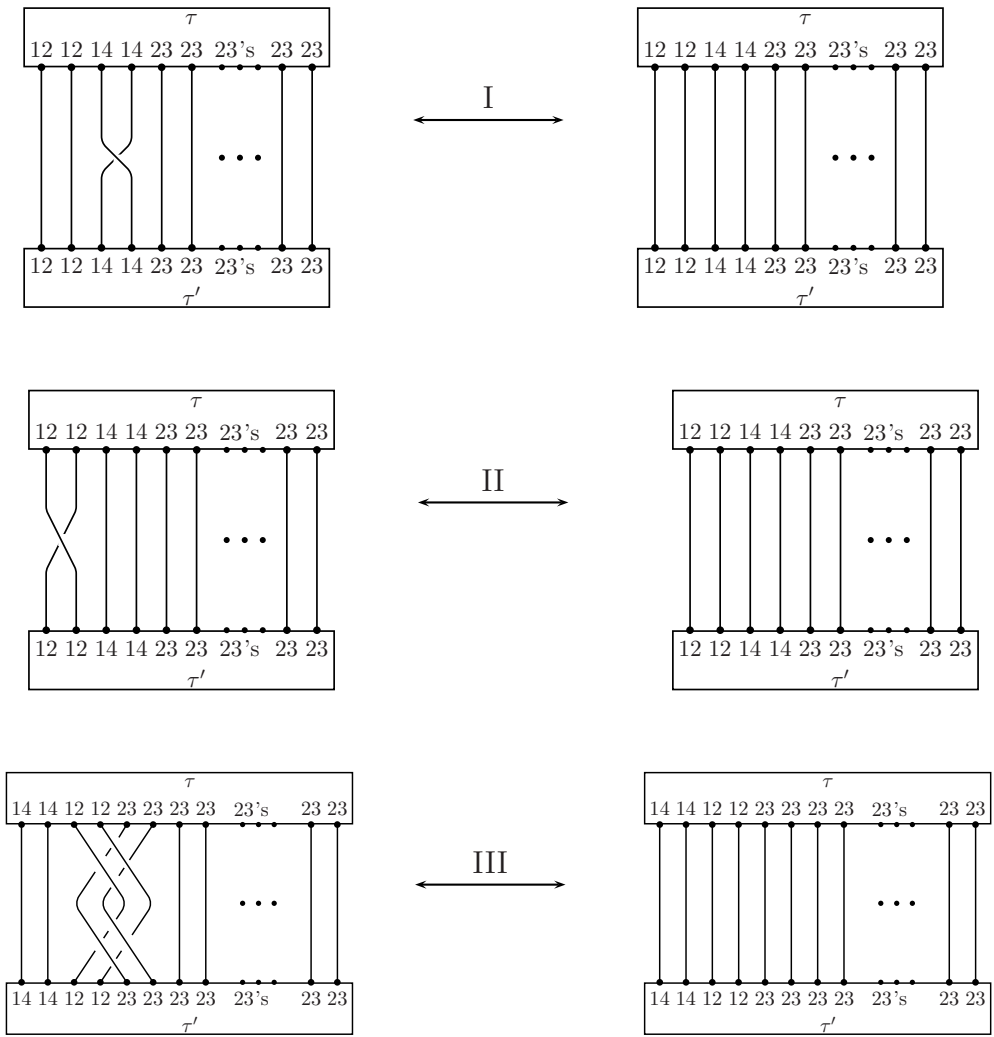

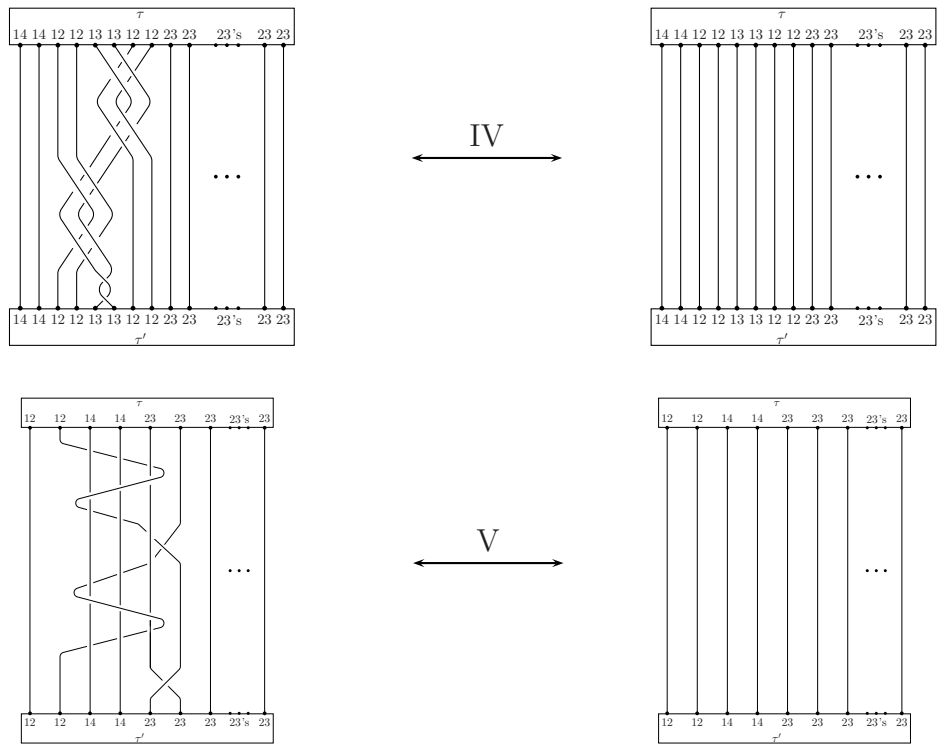
\begin{figure}[ht!]\vspace{-0.5cm}
\centerline{\noindent\hglue 1.5cm  \mbox{
      \subfigure{\scalebox{0.5}{ 
 \begin{pspicture}(1.5,-5)(6,7)
 \psdots(-1,4)(-.5,4)(0,4)(.5,4)(1,4)(1.5,4)(2,4)(2.5,4)(3,4)(3.5,4)(5.5,4)(6,4)
        (-1,-3)(-.5,-3)(0,-3)(.5,-3)(1,-3)(1.5,-3)(2,-3)(2.5,-3)(3,-3)(3.5,-3)(5.5,-3)(6,-3)
   \psline[linearc=.25](0,4)(0,1)(1,-.5)(.8,-.8)
   \psline[linearc=.25](0.5,4)(0.5,1)(1.5,-.5)(1.1,-1.1)
   \psline(3,4)(3,-3)
   \psline(3.5,4)(3.5,-3)
   \psline(5.5,4)(5.5,-3)
   \psline(6,4)(6,-3)
  \psline[linearc=.25](1,4)(2,2.5)(1.8,2.2)
  \psline[linearc=.25](1.5,4)(2.5,2.5)(2.1,1.9)
  \psline[linearc=.25](1.4,3.1)(1,2.5)(2,1)(2,-3)
  \psline[linearc=.25](1.7,2.8)(1.5,2.5)(2.5,1)(2.5,-3)
      \psline(1.6,3.4)(1.7,3.55)
   \psline(1.85,3.755)(2,4)
       \psline(1.85,3.025)(1.95,3.175)
    \psline(2.05,3.325)(2.5,4)
       \psline(1.65,1.975)(1.56,1.84)
    \psline(1.4,1.6)(.85,.755)
       \psline(1.9,1.6)(1.8,1.45)
    \psline(1.7,1.3)(1.5,1)
           \psframe(-1.3,-4)(6.3,-3)
           \psframe(-1.3,4)(6.3,5)
 \rput(-1,4.3){$14$}
 \rput(-.5,4.3){$14$}
 \rput(0,4.3){$12$}
 \rput(.5,4.3){$12$}
 \rput(1,4.3){$13$}
 \rput(1.5,4.3){$13$}
 \rput(2,4.3){$12$}
 \rput(2.5,4.3){$12$}
 \rput(3,4.3){$23$}
 \rput(3.5,4.3){$23$}
 \rput(5.5,4.3){$23$}
 \rput(6,4.3){$23$}
                   \rput(4.55,4.3){$23$'s}
 \rput(-1,-3.3){$14$}
 \rput(-.5,-3.3){$14$}
 \rput(0,-3.3){$12$}
 \rput(.5,-3.3){$12$}
 \rput(1,-3.3){$13$}
 \rput(1.5,-3.3){$13$}
 \rput(2,-3.3){$12$}
 \rput(2.5,-3.3){$12$}
 \rput(3,-3.3){$23$}
 \rput(3.5,-3.3){$23$}
 \rput(5.5,-3.3){$23$}
 \rput(6,-3.3){$23$}
                   \rput(4.55,-3.3){$23$'s}
  \psdots[dotsize=0.1](4.2,4)(4.5,4)(4.8,4)(4.2,-3)(4.5,-3)(4.8,-3)
                      (4.2,.5)(4.5,.5)(4.8,.5)
     \rput(2.2,4.8){$\tau$}
     \rput(2.2,-3.8){$\tau'$}
  \psline[linearc=.25](.7,-.2)(.5,-.5)(1.5,-2)(1.3,-2.2)
   \psline[linearc=.25](.4,.1)(0,-.5)(1,-2)(1.5,-2.5)(1.3,-2.7)
   \psline(.6,.4)(.7,.55)
      \psline(.85,.025)(.95,.175)
    \psline(1.05,.325)(1.5,1)
      \psline(.65,-1.025)(.56,-1.16)
    \psline[linearc=.25](.4,-1.4)(0,-2)(0,-3)
      \psline(.9,-1.4)(.8,-1.55)
    \psline[linearc=.25](.7,-1.7)(.5,-2)(.5,-3)
 \psline[linearc=.25](1.2,-2.3)(1,-2.5)(1.5,-3)
   \psline(1.2,-2.8)(1,-3)
  \psline(-.5,4)(-.5,-3)
  \psline(-1,4)(-1,-3)
\end{pspicture}
}
}
\quad
\begin{pspicture}(4,5)
\psline{<->}(1,3)(3,3)
\rput(2,3.3){\large{IV}}
\end{pspicture}
\quad
\subfigure{\scalebox{0.5}{
\begin{pspicture}(-.5,-5)(6.5,7)
\psdots(-1,4)(-.5,4)(0,4)(.5,4)(1,4)(1.5,4)(2,4)(2.5,4)(3,4)(3.5,4)(5.5,4)(6,4)
        (-1,-3)(-.5,-3)(0,-3)(.5,-3)(1,-3)(1.5,-3)(2,-3)(2.5,-3)(3,-3)(3.5,-3)
         (5.5,-3)(6,-3)
 \psline(-1,4)(-1,-3)
 \psline(-.5,4)(-.5,-3)
 \psline(0,4)(0,-3) 
 \psline(.5,4)(.5,-3)
 \psline(1,4)(1,-3)
\psline(1.5,4)(1.5,-3)
\psline(2,4)(2,-3)
\psline(2.5,4)(2.5,-3)
\psline(3,4)(3,-3)
\psline(3.5,4)(3.5,-3)
\psline(5.5,4)(5.5,-3)
\psline(6,4)(6,-3)
          \psframe(-1.3,-4)(6.3,-3)
           \psframe(-1.3,4)(6.3,5)
 \rput(-1,4.3){$14$}
 \rput(-.5,4.3){$14$}
 \rput(0,4.3){$12$}
 \rput(.5,4.3){$12$}
 \rput(1,4.3){$13$}
 \rput(1.5,4.3){$13$}
 \rput(2,4.3){$12$}
 \rput(2.5,4.3){$12$}
 \rput(3,4.3){$23$}
 \rput(3.5,4.3){$23$}
 \rput(5.5,4.3){$23$}
 \rput(6,4.3){$23$}
                   \rput(4.55,4.3){$23$'s}
 \rput(-1,-3.3){$14$}
 \rput(-.5,-3.3){$14$}
 \rput(0,-3.3){$12$}
 \rput(.5,-3.3){$12$}
 \rput(1,-3.3){$13$}
 \rput(1.5,-3.3){$13$}
 \rput(2,-3.3){$12$}
 \rput(2.5,-3.3){$12$}
 \rput(3,-3.3){$23$}
 \rput(3.5,-3.3){$23$}
 \rput(5.5,-3.3){$23$}
 \rput(6,-3.3){$23$}
                   \rput(4.55,-3.3){$23$'s}
  \psdots[dotsize=0.1](4.2,4)(4.5,4)(4.8,4)(4.2,-3)(4.5,-3)(4.8,-3)
                      (4.2,.5)(4.5,.5)(4.8,.5)
     \rput(2.2,4.8){$\tau$}
     \rput(2.2,-3.8){$\tau'$}
\end{pspicture} 
}
}
}}
\vspace{-2cm}
\centerline{\noindent\hglue 1.3cm\mbox{
      \subfigure{\scalebox{0.4}{ 
\begin{pspicture}(1.5,0)(6,10)
  \psline[linearc=.15](5,9.5)(5,6.3)(4,5)(1.1,4)(4.9,3)(1,2)(1,1.5)(1,0)
 \psframe[linecolor=white,fillstyle=solid,fillcolor=white](4.3,5.47)(4.5,5.61)
\psline[linearc=.15](4,3.1)(4,1.5)(5,.5)(5,0)
\psframe[linecolor=white,fillstyle=solid,fillcolor=white](4.4,.9)(4.6,1.1)
   \psline[linearc=.15](1,9.5)(1,9)(4.9,8)(1.1,7)(3.6,6.3)(5,5)(5,1.5)(4,.5)(4,0)
\psframe[linecolor=white,fillstyle=solid,fillcolor=white](1.9,8.65)(2.1,8.85)
\psframe[linecolor=white,fillstyle=solid,fillcolor=white](2.9,8.4)(3.1,8.6)
\psframe[linecolor=white,fillstyle=solid,fillcolor=white](3.9,8.15)(4.1,8.35)
\psframe[linecolor=white,fillstyle=solid,fillcolor=white](1.9,6.7)(2.1,6.9)
\psframe[linecolor=white,fillstyle=solid,fillcolor=white](2.9,6.4)(3.1,6.6)
\psframe[linecolor=white,fillstyle=solid,fillcolor=white](3.85,5.8)(4.15,6.05)
\psframe[linecolor=white,fillstyle=solid,fillcolor=white](1.9,4.2)(2.1,4.4)
\psframe[linecolor=white,fillstyle=solid,fillcolor=white](2.9,4.5)(3.1,4.7)
\psframe[linecolor=white,fillstyle=solid,fillcolor=white](3.9,4.9)(4.15,5.3)
\psframe[linecolor=white,fillstyle=solid,fillcolor=white](1.9,2.1)(2.1,2.4)
\psframe[linecolor=white,fillstyle=solid,fillcolor=white](2.9,2.4)(3.1,2.6)
\psframe[linecolor=white,fillstyle=solid,fillcolor=white](3.9,2.7)(4.1,2.9)
  \psline(2,9.5)(2,7.35)   \psline(2,7.12)(2,3.9)    \psline(2,3.65)(2,0) 
  \psline(3,9.5)(3,7.6)    \psline(3,7.4)(3,3.6)     \psline(3,3.4)(3,0)
  \psline(4,9.5)(4,7.9)    \psline(4,7.6)(4,3.35)    \psline(4,3.1)(4,1.5)
  \psline(0,0)(0,9.5)      \psline(6,0)(6,9.5)       \psline(7.8,9.5)(7.8,0)
 \psdots(0,9.5)(1,9.5)(2,9.5)(3,9.5)(4,9.5)(5,9.5)(6,9.5)(7.8,9.5)
        (0,0)(1,0)(2,0)(3,0)(4,0)(5,0)(6,0)(7.8,0)
                 \psdots[dotsize=.1](6.65,9.5)(6.95,9.5)(7.25,9.5)
                   \rput(6.95,9.8){$23$'s}
                \psdots[dotsize=.1](6.65,0)(6.95,0)(7.25,0)
                   \rput(6.95,-.3){$23$'s}
                \psdots[dotsize=.1](6.65,4.75)(6.95,4.75)(7.25,4.75)
       \psframe(-.3,9.5)(8.1,10.5)
       \psframe(-.3,0)(8.1,-1)
      \rput(3.65,10.3){$\tau$}
      \rput(3.65,-.8){$\tau'$}
 \rput(0,9.8){$12$}
 \rput(1,9.8){$12$}
 \rput(2,9.8){$14$}  
 \rput(3,9.8){$14$}  
 \rput(4,9.8){$23$}
 \rput(5,9.8){$23$}
 \rput(6,9.8){$23$}
 \rput(7.8,9.8){$23$}
 \rput(0,-.3){$12$}
 \rput(1,-.3){$12$}
 \rput(2,-.3){$14$}
 \rput(3,-.3){$14$}
 \rput(4,-.3){$23$}
 \rput(5,-.3){$23$}
 \rput(6,-.3){$23$}
 \rput(7.8,-.3){$23$}
\end{pspicture}
}
}
\quad
\begin{pspicture}(5,6)
\psline{<->}(2,2)(4,2)
\rput(3,2.3){\large{V}}
\end{pspicture}
 \subfigure{\scalebox{0.4}{ 
\begin{pspicture}(-.5,0)(10,10)
  \psline(1,9.5)(1,0)      \psline(4,9.5)(4,0)       \psline(5,9.5)(5,0)
   \psline(0,0)(0,9.5)      \psline(6,0)(6,9.5)       \psline(7.8,9.5)(7.8,0)
\psline(3,9.5)(3,0)
\psline(2,9.5)(2,0)
 \psdots(0,9.5)(1,9.5)(2,9.5)(3,9.5)(4,9.5)(5,9.5)(6,9.5)(7.8,9.5)
        (0,0)(1,0)(2,0)(3,0)(4,0)(5,0)(6,0)(7.8,0)
                 \psdots[dotsize=.1](6.65,9.5)(6.95,9.5)(7.25,9.5)
                   \rput(6.95,9.8){$23$'s}
                \psdots[dotsize=.1](6.65,0)(6.95,0)(7.25,0)
                   \rput(6.95,-.3){$23$'s}
                \psdots[dotsize=.1](6.65,4.75)(6.95,4.75)(7.25,4.75)
       \psframe(-.3,9.5)(8.1,10.5)
       \psframe(-.3,0)(8.1,-1)
      \rput(3.65,10.3){$\tau$}
      \rput(3.65,-.8){$\tau'$}
 \rput(0,9.8){$12$}
 \rput(1,9.8){$12$}
 \rput(2,9.8){$14$}  
 \rput(3,9.8){$14$}  
 \rput(4,9.8){$23$}
 \rput(5,9.8){$23$}
 \rput(6,9.8){$23$}
 \rput(7.8,9.8){$23$}
 \rput(0,-.3){$12$}
 \rput(1,-.3){$12$}
 \rput(2,-.3){$14$}
 \rput(3,-.3){$14$}
 \rput(4,-.3){$23$}
 \rput(5,-.3){$23$}
 \rput(6,-.3){$23$}
 \rput(7.8,-.3){$23$}
\end{pspicture}
}
}
}}
\vspace{0.5cm}  
\caption{The last two non-local moves}
  \label{fig:movesII}
\end{figure}

\begin{thm}\label{thm:movesI}
 Two bi-tricolored links represent the  same $3$--manifold if and 
only if they can be related 
using a finite number of moves from the following list: 
\begin{itemize}
\item bi-tricolored Reidemeister moves
\item moves $\CM$, $\CP$, (see Figure~\ref{fig:localmoves})
\item moves I,\ldots,V shown in Figures~\ref{fig:movesI}, \ref{fig:movesII}.
\end{itemize}  
\end{thm}
\begin{rem}
  Notice the similarity between moves II, III, and IV and the ``homonymous'' moves 
in~\cite{Pier1991}. This similarity will be  exploited in the proof of 
this and the next theorems.     
\end{rem}

\proof
  Since the braids on the left side of moves I-V belong to the kernel of the 
lifting
homeomorphism, these  indeed do not change the represented manifold.\n
  The proof of the reverse direction is completely analogous to Piergallini's
 proof of the 
main theorem  of~\cite{Pier1991}.
We just comment on how each step of that
proof goes through in the present situation:
\begin{itemize}
\item  Each bi-tircolored link has a normalized diagram that is a plat diagram whose
top and bottom are colored by the standard $4$-fold covering. This is easily seen
by observing that the proof of the transitivity of the braid action 
in~\cite{BernEdm1979} can be carried out by colored isotopy, much as Piergallini does
in the $3$-fold case.
\item The first step (Heegaard stabilization) can be realized in exactly the 
same way.
\item For the second step (Heegaard equivalence), in order to get braids that 
lift to the Suzuki
 generators
one just has to add two trivial strands colored by $(1,4)$ to the braids 
Piergallini uses.
\item For the third step, we have to show that one can add
     each element of $\ker\bar{\lambda}$ to the top (and bottom)
     of normalized diagrams using the moves. This will be
     accomplished if we show that one can insert in a normalized
     diagram, in any position with the right colors, each of the
     normal generators of $\ker\bar{\lambda}$ given in
     Theorem~\ref{thm:kernel}.\footnote{This proof was corrected 20 January 2004.  I am grateful to Riccardo Piergallini for pointing out the mistake.}\n
Braids $\Gb_0$ and $\Gb_2$ can be obviously added using moves I and II, respectively. 
Also by slight 
modification of
Piegallini's proof (by just adding two trivial strands colored by $(1,4)$), $B_1$ 
and $D_1$ can 
be added using
moves $\mathcal{M}$, I-IV.\n
Now move V transforms $\delta_4$ to $\Gb_4$
 and therefore using move V  we can transform the remaining
  normal generators of $\ker\bar{\lambda}$ to $B$ or $D$ (it is easily checked that 
the colors are 
right).
     Thus any element of
$\ker\bar{\lambda}$ can be added on the top and the bottom of a normalized 
diagram using
 the given moves.\qed
\end{itemize}

The following combination of move $\CM$ and isotopy will be used throughout what 
follows, so we give it a name.
\begin{lem}\label{newref} 
The ``circumcision''  move shown in Figure~\ref{fig:circum} is a 
consequence of move
$\mathcal{M}$.
\begin{figure}[ht!]\small
\centering
\begin{pspicture}(0,-.7)(3.5,2.7)
 \rput(1,2.8){$ij$} \rput(2,2.8){$ij$} \rput(1,-.7){$ik$} \rput(2,-.7){$ik$}
   \psellipse(1.5,1)(.7,.4)
    \psframe[linecolor=white,fillstyle=solid,fillcolor=white](.9,1.1)(1.1,1.4)
    \psframe[linecolor=white,fillstyle=solid,fillcolor=white](1.9,1.1)(2.1,1.4)
 \psline(1,2.5)(1,.85) \psline(1,.65)(1,-.5)
   \psline(2,2.5)(2,.85) \psline(2,.65)(2,-.5)
  \rput(.5,1){$jk$}
\end{pspicture}
\begin{pspicture}(.2,-.7)(2.9,2)
  \psline{<->}(.7,1.5)(2.3,1.5)
\end{pspicture}
\begin{pspicture}(0,-.7)(3.5,2.7)
 \rput(1,2.8){$ij$} \rput(2,2.8){$ij$} \rput(1,-.7){$ik$}
\rput(2,-.7){$ik$}
 \psline(1,2.5)(1,2) \psline(2,2.5)(2,2)
      \psarc(1.5,2){.5}{180}{360}
      \psarc(1.5,0){.5}{0}{180}
 \psline(1,-.5)(1,0) \psline(2,-.5)(2,0)
\end{pspicture}
\caption{The circumcision move}\label{fig:circum}
\end{figure}
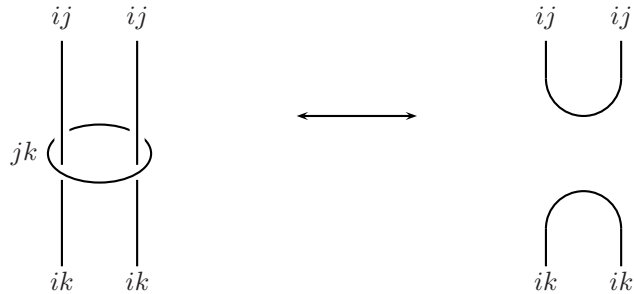
\end{lem}

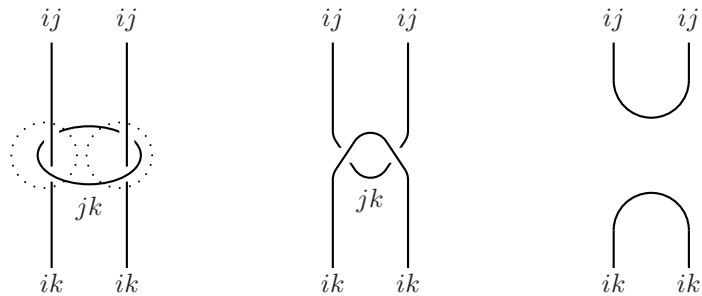
\begin{figure}[ht!]\small\vspace{1cm}
\cl{\begin{pspicture}(0,-1.5)(2.5,2.7)
 \rput(1,2.8){$ij$} \rput(2,2.8){$ij$} \rput(1,-.7){$ik$} \rput(2,-.7){$ik$}
   \psellipse(1.5,1)(.7,.4)
    \psframe[linecolor=white,fillstyle=solid,fillcolor=white](.9,1.1)(1.1,1.4)
    \psframe[linecolor=white,fillstyle=solid,fillcolor=white](1.9,1.1)(2.1,1.4)
 \psline(1,2.5)(1,.85) \psline(1,.65)(1,-.5)
   \psline(2,2.5)(2,.85) \psline(2,.65)(2,-.5)
  \rput(1.5,.3){$jk$}
 \pscircle[linestyle=dotted](.9,1){.45}
 \pscircle[linestyle=dotted](1.9,1){.45}
\end{pspicture}
\hspace{1cm}
\begin{pspicture}(0,-1.5)(2.5,2.7)
 \rput(1,2.8){$ij$} \rput(2,2.8){$ij$} \rput(1,-.7){$ik$} \rput(2,-.7){$ik$}
  \psline[linearc=.25](1,2.5)(1,1.25)(1.5,.5)(2,1.25)(2,2.5)
\psframe[linecolor=white,fillstyle=solid,fillcolor=white](1,.9)(2,1.1)
  \psline[linearc=.25](1,-.5)(1,.75)(1.5,1.5)(2,.75)(2,-.5)
  \rput(1.5,.4){$jk$}
\end{pspicture}
 \hspace{1cm}
\begin{pspicture}(0,-1.5)(2.5,2.7)
 \rput(1,2.8){$ij$} \rput(2,2.8){$ij$} \rput(1,-.7){$ik$}
\rput(2,-.7){$ik$}
 \psline(1,2.5)(1,2) \psline(2,2.5)(2,2)
      \psarc(1.5,2){.5}{180}{360}
      \psarc(1.5,0){.5}{0}{180}
 \psline(1,-.5)(1,0) \psline(2,-.5)(2,0)
\end{pspicture}}
\vspace{-0.3cm}
\caption{Proof of lemma \ref{newref}}\label{newfig}
\end{figure}

  \begin{proof}
 Perform $\mathcal{M}$ moves inside the dotted circles and then
 isotope as in figure \ref{newfig}.\end{proof}

\begin{thm}
\label{thm:main}
  Moves $\CM$, $\CP$ and stabilization suffice. Actually one needs only five sheets. 
\end{thm}

\begin{proof}
  It suffices to show how the non-local moves I,\ldots,V can be realized. 
For moves I and V this shown in Figures~\ref{fig:Iproof} and \ref{fig:Vproof} 
respectively.\n
For moves II, III and IV we use the similarity with the homonymous moves 
of ~\cite{Pier1991}: the proofs for these are very similar to Piergallini's  
proofs in~\cite{Pier1995}. The new element is that one can get a circumcising 
fourth sheet ``for free'' and thus five sheets suffice. 
This is done in detail for move II in Figure~\ref{fig:IIproof}
 leaving the remaining cases to the reader.    
\end{proof}

\begin{figure}[htp]
  \begin{center}
    \mbox{
      \subfigure[Start by adding a trivial sheet with monodromy {$(4,5)$}]
{\scalebox{0.7}{
\begin{pspicture}(.5,0)(6,5)
 \psdots(0,4)(.5,4)(1,4)(1.5,4)(2,4)(2.5,4)(4,4)(4.5,4)
        (0,1)(.5,1)(1,1)(1.5,1)(2,1)(2.5,1)(4,1)(4.5,1)
   \psline(0,4)(0,1)
   \psline(0.5,4)(0.5,1)
   \psline(2,4)(2,1)
   \psline(2.5,4)(2.5,1)
   \psline(4,4)(4,1)
   \psline(4.5,4)(4.5,1)
 \psline[linearc=.25](1.5,4)(1.5,2.75)(1,2.25)(1,1)
\psframe[linecolor=white,fillstyle=solid,fillcolor=white](1.2,2.4)(1.38,2.55)
\psline[linearc=.25](1,4)(1,2.75)(1.5,2.25)(1.5,1)
           \psframe(-.3,0)(4.8,1)
           \psframe(-.3,4)(4.8,5)
    \rput(0,4.3){$12$}
 \rput(0.5,4.3){$12$}
 \rput(1,4.3){$14$}
 \rput(1.5,4.3){$14$}
 \rput(2,4.3){$23$}
 \rput(2.5,4.3){$23$}
 \rput(4,4.3){$23$}
 \rput(4.5,4.3){$23$}
                   \rput(3.2,4.3){$23$'s}
 \rput(0,.7){$12$}
 \rput(0.5,.7){$12$}
 \rput(1,.7){$14$}
 \rput(1.5,.7){$14$}
 \rput(2,.7){$23$}
 \rput(2.5,.7){$23$}
 \rput(4,.7){$23$}
 \rput(4.5,.7){$23$}
                   \rput(3.2,.7){$23$'s}
  \psdots[dotsize=0.1](3,4)(3.3,4)(3.6,4)(3,1)(3.3,1)(3.6,1)
                      (3,2.5)(3.3,2.5)(3.6,2.5)
     \rput(2.5,4.8){$\tau$}
     \rput(2.5,.2){$\tau'$}
    \pscircle(6,2.5){.5}
  \rput(6.8,2.5){$45$}
\end{pspicture} 
}
}
 \qquad \quad
      \subfigure[flip it over the whole diagram\ldots]{\scalebox{0.7}{
\begin{pspicture}(-.5,0)(6,5)
 \psdots(0,4)(.5,4)(1,4)(1.5,4)(2,4)(2.5,4)(4,4)(4.5,4)
        (0,1)(.5,1)(1,1)(1.5,1)(2,1)(2.5,1)(4,1)(4.5,1)
   \psline(0,4)(0,1)
   \psline(0.5,4)(0.5,1)
   \psline(2,4)(2,1)
   \psline(2.5,4)(2.5,1)
   \psline(4,4)(4,1)
   \psline(4.5,4)(4.5,1)
   \psline(4.5,4)(4.5,1)
 \psline[linearc=.25](1.5,4)(1.5,2.75)(1,2.25)(1,1)
\psframe[linecolor=white,fillstyle=solid,fillcolor=white](1.2,2.4)(1.38,2.55)
\psline[linearc=.25](1,4)(1,2.75)(1.5,2.25)(1.5,1)
           \psframe(-.3,0)(4.8,1)
           \psframe(-.3,4)(4.8,5)
    \rput(0,4.3){$12$}
 \rput(0.5,4.3){$12$}
 \rput(1,4.3){$15$}
 \rput(1.5,4.3){$15$}
 \rput(2,4.3){$23$}
 \rput(2.5,4.3){$23$}
 \rput(4,4.3){$23$}
 \rput(4.5,4.3){$23$}
                   \rput(3.2,4.3){$23$'s}
 \rput(0,.7){$12$}
 \rput(0.5,.7){$12$}
 \rput(1,.7){$14$}
 \rput(1.5,.7){$14$}
 \rput(2,.7){$23$}
 \rput(2.5,.7){$23$}
 \rput(4,.7){$23$}
 \rput(4.5,.7){$23$}
                   \rput(3.2,.7){$23$'s}
  \psdots[dotsize=0.1](3,4)(3.3,4)(3.6,4)(3,1)(3.3,1)(3.6,1)
                      (3,2.5)(3.3,2.5)(3.6,2.5)
     \rput(2.5,4.8){$\tau$}
     \rput(2.5,.2){$\tau'$}
\psframe[linecolor=white,fillstyle=solid,fillcolor=white](-1,3.05)(5.5,3.35)
  \psellipse(2.3,3.4)(2.9,.27)
  \psframe[linecolor=white,fillstyle=solid,fillcolor=white](-.1,3.4)(.1,3.8)
\psframe[linecolor=white,fillstyle=solid,fillcolor=white](.4,3.4)(.6,3.8)
\psframe[linecolor=white,fillstyle=solid,fillcolor=white](.9,3.4)(1.1,3.8)
\psframe[linecolor=white,fillstyle=solid,fillcolor=white](1.4,3.4)(1.6,3.8)
\psframe[linecolor=white,fillstyle=solid,fillcolor=white](1.9,3.4)(2.1,3.8)
\psframe[linecolor=white,fillstyle=solid,fillcolor=white](2.4,3.4)(2.6,3.8)
\psframe[linecolor=white,fillstyle=solid,fillcolor=white](3.9,3.4)(4.1,3.8)
\psframe[linecolor=white,fillstyle=solid,fillcolor=white](4.4,3.4)(4.6,3.8)
\psline(0,4)(0,3.38)
\psline(.5,4)(.5,3.36)
\psline(1,4)(1,3.3)
\psline(1.5,4)(1.5,3.3)
\psline(2,4)(2,3.32)
\psline(2.5,4)(2.5,3.32)
\psline(4,4)(4,3.36)
\psline(4.5,4)(4.5,3.38)
\rput(5.5,3.4){$45$}
\end{pspicture}
}%
}%
}%
\vspace{.2cm}
    \mbox{
      \subfigure[and apply a sequence of {$\CP$} moves\ldots]{\scalebox{0.7}{
\begin{pspicture}(.5,0)(6,5)
 \psdots(0,4)(.5,4)(1,4)(1.5,4)(2,4)(2.5,4)(4,4)(4.5,4)
        (0,1)(.5,1)(1,1)(1.5,1)(2,1)(2.5,1)(4,1)(4.5,1)
   \psline(0,4)(0,1)
   \psline(0.5,4)(0.5,1)
   \psline(2,4)(2,1)
   \psline(2.5,4)(2.5,1)
   \psline(4,4)(4,1)
   \psline(4.5,4)(4.5,1)
 \psline[linearc=.25](1.5,4)(1.5,2.75)(1,2.25)(1,1)
\psframe[linecolor=white,fillstyle=solid,fillcolor=white](1.2,2.4)(1.38,2.55)
\psline[linearc=.25](1,4)(1,2.75)(1.5,2.25)(1.5,1)
           \psframe(-.3,0)(4.8,1)
           \psframe(-.3,4)(4.8,5)
    \rput(0,4.3){$12$}
 \rput(0.5,4.3){$12$}
 \rput(1,4.3){$15$}
 \rput(1.5,4.3){$15$}
 \rput(2,4.3){$23$}
 \rput(2.5,4.3){$23$}
 \rput(4,4.3){$23$}
 \rput(4.5,4.3){$23$}
                   \rput(3.2,4.3){$23$'s}
 \rput(0,.7){$12$}
 \rput(0.5,.7){$12$}
 \rput(1,.7){$14$}
 \rput(1.5,.7){$14$}
 \rput(2,.7){$23$}
 \rput(2.5,.7){$23$}
 \rput(4,.7){$23$}
 \rput(4.5,.7){$23$}
                   \rput(3.2,.7){$23$'s}
  \psdots[dotsize=0.1](3,4)(3.3,4)(3.6,4)(3,1)(3.3,1)(3.6,1)
                      (3,2.5)(3.3,2.5)(3.6,2.5)
     \rput(2.5,4.8){$\tau$}
     \rput(2.5,.2){$\tau'$}
\psframe[linecolor=white,fillstyle=solid,fillcolor=white](.9,3.1)(1.6,3.4)
    \psellipse(1.25,3.4)(.5,.3)
 \psframe[linecolor=white,fillstyle=solid,fillcolor=white](.9,3.5)(1.1,3.8)
\psframe[linecolor=white,fillstyle=solid,fillcolor=white](1.4,3.5)(1.6,3.8)
\psline(1,4)(1,3.3)
\psline(1.5,4)(1.5,3.3)
\end{pspicture}
}
}
 \qquad \quad
      \subfigure[circumcise \ldots]{\scalebox{0.7}{
\begin{pspicture}(-.5,0)(5,5)
 \psdots(0,4)(.5,4)(1,4)(1.5,4)(2,4)(2.5,4)(4,4)(4.5,4)
        (0,1)(.5,1)(1,1)(1.5,1)(2,1)(2.5,1)(4,1)(4.5,1)
   \psline(0,4)(0,1)
   \psline(0.5,4)(0.5,1)
   \psline(2,4)(2,1)
   \psline(2.5,4)(2.5,1)
   \psline(4,4)(4,1)
   \psline(4.5,4)(4.5,1)
 \psline[linearc=.25](1.01,3)(1.25,3.44)(1.5,3)(1.5,2.75)(1,2.25)(1,1)
\psframe[linecolor=white,fillstyle=solid,fillcolor=white](1.2,2.4)(1.38,2.55)
\psline[linearc=.25](1.491,3)(1.25,3.44)(1,3)(1,2.75)(1.5,2.25)(1.5,1)
\psline[linearc=.25](1,4)(1,3.5)(1.25,3.1)(1.5,3.5)(1.5,4)
           \psframe(-.3,0)(4.8,1)
           \psframe(-.3,4)(4.8,5)
    \rput(0,4.3){$12$}
 \rput(0.5,4.3){$12$}
 \rput(1,4.3){$15$}
 \rput(1.5,4.3){$15$}
 \rput(2,4.3){$23$}
 \rput(2.5,4.3){$23$}
 \rput(4,4.3){$23$}
 \rput(4.5,4.3){$23$}
                   \rput(3.2,4.3){$23$'s}
 \rput(0,.7){$12$}
 \rput(0.5,.7){$12$}
 \rput(1,.7){$14$}
 \rput(1.5,.7){$14$}
 \rput(2,.7){$23$}
 \rput(2.5,.7){$23$}
 \rput(4,.7){$23$}
 \rput(4.5,.7){$23$}
                   \rput(3.2,.7){$23$'s}
  \psdots[dotsize=0.1](3,4)(3.3,4)(3.6,4)(3,1)(3.3,1)(3.6,1)
                      (3,2.5)(3.3,2.5)(3.6,2.5)
     \rput(2.5,4.8){$\tau$}
     \rput(2.5,.2){$\tau'$}
\end{pspicture}
}%
}%
}%
\vspace{.2cm}
\mbox{\subfigure[and isotope to get rid of the crossing.]{\scalebox{0.7}{
\begin{pspicture}(-.5,0)(5,5)
 \psdots(0,4)(.5,4)(1,4)(1.5,4)(2,4)(2.5,4)(4,4)(4.5,4)
        (0,1)(.5,1)(1,1)(1.5,1)(2,1)(2.5,1)(4,1)(4.5,1)
   \psline(0,4)(0,1)
   \psline(0.5,4)(0.5,1)
   \psline(2,4)(2,1)
   \psline(2.5,4)(2.5,1)
   \psline(4,4)(4,1)
   \psline(4.5,4)(4.5,1)
 \psline[linearc=.25](1,1)(1,3)(1.25,3.4)(1.5,3)(1.5,1)
\psline[linearc=.25](1,4)(1,3.5)(1.25,3.1)(1.5,3.5)(1.5,4)
           \psframe(-.3,0)(4.8,1)
           \psframe(-.3,4)(4.8,5)
    \rput(0,4.3){$12$}
 \rput(0.5,4.3){$12$}
 \rput(1,4.3){$15$}
 \rput(1.5,4.3){$15$}
 \rput(2,4.3){$23$}
 \rput(2.5,4.3){$23$}
 \rput(4,4.3){$23$}
 \rput(4.5,4.3){$23$}
                   \rput(3.2,4.3){$23$'s}
 \rput(0,.7){$12$}
 \rput(0.5,.7){$12$}
 \rput(1,.7){$14$}
 \rput(1.5,.7){$14$}
 \rput(2,.7){$23$}
 \rput(2.5,.7){$23$}
 \rput(4,.7){$23$}
 \rput(4.5,.7){$23$}
                   \rput(3.2,.7){$23$'s}
  \psdots[dotsize=0.1](3,4)(3.3,4)(3.6,4)(3,1)(3.3,1)(3.6,1)
                      (3,2.5)(3.3,2.5)(3.6,2.5)
     \rput(2.5,4.8){$\tau$}
     \rput(2.5,.2){$\tau'$}
\end{pspicture}
}}
\qquad \qquad \quad 
\subfigure[Now reverse the process starting at (d) to complete move I.]{\scalebox{0.7}{
\begin{pspicture}(.5,0)(6,5)
  \psdots(0,4)(.5,4)(1,4)(1.5,4)(2,4)(2.5,4)(4,4)(4.5,4)
        (0,1)(.5,1)(1,1)(1.5,1)(2,1)(2.5,1)(4,1)(4.5,1)
   \psline(0,4)(0,1)
   \psline(0.5,4)(0.5,1)
   \psline(2,4)(2,1)
   \psline(2.5,4)(2.5,1)
   \psline(4,4)(4,1)
   \psline(4.5,4)(4.5,1)
 \psline(1.5,4)(1.5,1)
\psline(1,4)(1,1)
           \psframe(-.3,0)(4.8,1)
           \psframe(-.3,4)(4.8,5)
    \rput(0,4.3){$12$}
 \rput(0.5,4.3){$12$}
 \rput(1,4.3){$14$}
 \rput(1.5,4.3){$14$}
 \rput(2,4.3){$23$}
 \rput(2.5,4.3){$23$}
 \rput(4,4.3){$23$}
 \rput(4.5,4.3){$23$}
                   \rput(3.2,4.3){$23$'s}
 \rput(0,.7){$12$}
 \rput(0.5,.7){$12$}
 \rput(1,.7){$14$}
 \rput(1.5,.7){$14$}
 \rput(2,.7){$23$}
 \rput(2.5,.7){$23$}
 \rput(4,.7){$23$}
 \rput(4.5,.7){$23$}
                   \rput(3.2,.7){$23$'s}
  \psdots[dotsize=0.1](3,4)(3.3,4)(3.6,4)(3,1)(3.3,1)(3.6,1)
                      (3,2.5)(3.3,2.5)(3.6,2.5)
     \rput(2.5,4.8){$\tau$}
     \rput(2.5,.2){$\tau'$}
\end{pspicture} 
      }}}
    \caption{Proof of Theorem~\ref{thm:main} for move I}
    \label{fig:Iproof}
  \end{center}
\end{figure}
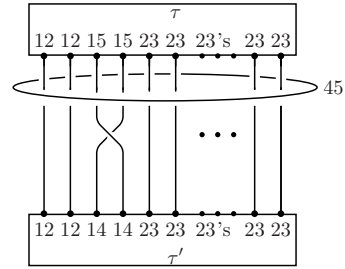
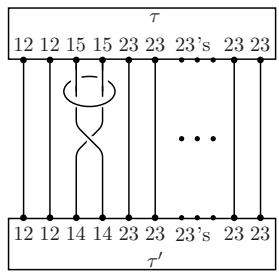
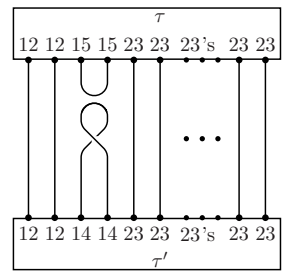
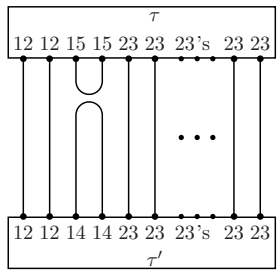
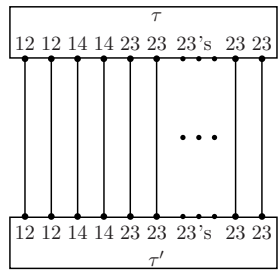

\begin{figure}[htp]\vspace{-10pt}
  \centering
  \mbox{
\subfigure
{ \scalebox{.5}{
\begin{pspicture}(-.5,-1)(10,10)
  \psline[linearc=.15](5,9.5)(5,6.3)(4,5)(1.1,4)(4.9,3)(1,2)(1,1.5)(1,0)
 \psframe[linecolor=white,fillstyle=solid,fillcolor=white](4.3,5.47)(4.5,5.61)
\psline[linearc=.15](4,3.1)(4,1.5)(5,.5)(5,0)
\psframe[linecolor=white,fillstyle=solid,fillcolor=white](4.4,.9)(4.6,1.1)
   \psline[linearc=.15](1,9.5)(1,9)(4.9,8)(1.1,7)(3.6,6.3)(5,5)(5,1.5)(4,.5)(4,0)
\psframe[linecolor=white,fillstyle=solid,fillcolor=white](1.9,8.65)(2.1,8.85)
\psframe[linecolor=white,fillstyle=solid,fillcolor=white](2.9,8.4)(3.1,8.6)
\psframe[linecolor=white,fillstyle=solid,fillcolor=white](3.9,8.15)(4.1,8.35)
\psframe[linecolor=white,fillstyle=solid,fillcolor=white](1.9,6.7)(2.1,6.9)
\psframe[linecolor=white,fillstyle=solid,fillcolor=white](2.9,6.4)(3.1,6.6)
\psframe[linecolor=white,fillstyle=solid,fillcolor=white](3.85,5.8)(4.15,6.05)
\psframe[linecolor=white,fillstyle=solid,fillcolor=white](1.9,4.2)(2.1,4.4)
\psframe[linecolor=white,fillstyle=solid,fillcolor=white](2.9,4.5)(3.1,4.7)
\psframe[linecolor=white,fillstyle=solid,fillcolor=white](3.9,4.9)(4.15,5.3)
\psframe[linecolor=white,fillstyle=solid,fillcolor=white](1.9,2.1)(2.1,2.4)
\psframe[linecolor=white,fillstyle=solid,fillcolor=white](2.9,2.4)(3.1,2.6)
\psframe[linecolor=white,fillstyle=solid,fillcolor=white](3.9,2.7)(4.1,2.9)
  \psline(2,9.5)(2,7.35)   \psline(2,7.12)(2,3.9)    \psline(2,3.65)(2,0) 
  \psline(3,9.5)(3,7.6)    \psline(3,7.4)(3,3.6)     \psline(3,3.4)(3,0)
  \psline(4,9.5)(4,7.9)    \psline(4,7.6)(4,3.35)    \psline(4,3.1)(4,1.5)
  \psline(0,0)(0,9.5)      \psline(6,0)(6,9.5)       \psline(7.8,9.5)(7.8,0)
 \psdots(0,9.5)(1,9.5)(2,9.5)(3,9.5)(4,9.5)(5,9.5)(6,9.5)(7.8,9.5)
        (0,0)(1,0)(2,0)(3,0)(4,0)(5,0)(6,0)(7.8,0)
                 \psdots[dotsize=.1](6.65,9.5)(6.95,9.5)(7.25,9.5)
                   \rput(6.95,9.8){$23$'s}
                \psdots[dotsize=.1](6.65,0)(6.95,0)(7.25,0)
                   \rput(6.95,-.3){$23$'s}
                \psdots[dotsize=.1](6.65,4.75)(6.95,4.75)(7.25,4.75)
       \psframe(-.3,9.5)(8.1,10.5)
       \psframe(-.3,0)(8.1,-1)
      \rput(3.65,10.3){$\tau$}
      \rput(3.65,-.8){$\tau'$}
 \rput(0,9.8){$12$}
 \rput(1,9.8){$12$}
 \rput(2,9.8){$15$}  
 \rput(3,9.8){$15$}  
 \rput(4,9.8){$23$}
 \rput(5,9.8){$23$}
 \rput(6,9.8){$23$}
 \rput(7.8,9.8){$23$}
 \rput(0,-.3){$12$}
 \rput(1,-.3){$12$}
 \rput(2,-.3){$14$}
 \rput(3,-.3){$14$}
 \rput(4,-.3){$23$}
 \rput(5,-.3){$23$}
 \rput(6,-.3){$23$}
 \rput(7.8,-.3){$23$}
\psframe[linecolor=white,fillstyle=solid,fillcolor=white](1.8,5)(3.2,6)    %
\psellipse(2.5,5.5)(.8,.3)                                                 %
\rput(1.35,5.5){$45$}                                                      
\psframe[linecolor=white,fillstyle=solid,fillcolor=white](1.9,5.6)(2.1,5.8)
\psframe[linecolor=white,fillstyle=solid,fillcolor=white](2.9,5.6)(3.1,5.8)
\psline(2,6.2)(2,5.39)       \psline(2,5.17)(2,4.4)                        %
\psline(3,6.2)(3,5.39)       \psline(3,5.17)(3,4.4)                        %
\end{pspicture}
}}%
 \qquad \qquad
\subfigure
{ \scalebox{0.5}{
\begin{pspicture}(-.5,-1)(10,10)
  \psline[linearc=.15](5,9.5)(5,6.3)(4,5)(1.1,4)(4.9,3)(1,2)(1,1.5)(1,0)
 \psframe[linecolor=white,fillstyle=solid,fillcolor=white](4.3,5.47)(4.5,5.61)
\psline[linearc=.15](4,3.1)(4,1.5)(5,.5)(5,0)
\psframe[linecolor=white,fillstyle=solid,fillcolor=white](4.4,.9)(4.6,1.1)
   \psline[linearc=.15](1,9.5)(1,9)(4.9,8)(1.1,7)(3.6,6.3)(5,5)(5,1.5)(4,.5)(4,0)
\psframe[linecolor=white,fillstyle=solid,fillcolor=white](1.9,8.65)(2.1,8.85)
\psframe[linecolor=white,fillstyle=solid,fillcolor=white](2.9,8.4)(3.1,8.6)
\psframe[linecolor=white,fillstyle=solid,fillcolor=white](3.9,8.15)(4.1,8.35)
\psframe[linecolor=white,fillstyle=solid,fillcolor=white](1.9,6.7)(2.1,6.9)
\psframe[linecolor=white,fillstyle=solid,fillcolor=white](2.9,6.4)(3.1,6.6)
\psframe[linecolor=white,fillstyle=solid,fillcolor=white](3.85,5.8)(4.15,6.05)
\psframe[linecolor=white,fillstyle=solid,fillcolor=white](1.9,4.2)(2.1,4.4)
\psframe[linecolor=white,fillstyle=solid,fillcolor=white](2.9,4.5)(3.1,4.7)
\psframe[linecolor=white,fillstyle=solid,fillcolor=white](3.9,4.9)(4.15,5.3)
\psframe[linecolor=white,fillstyle=solid,fillcolor=white](1.9,2.1)(2.1,2.4)
\psframe[linecolor=white,fillstyle=solid,fillcolor=white](2.9,2.4)(3.1,2.6)
\psframe[linecolor=white,fillstyle=solid,fillcolor=white](3.9,2.7)(4.1,2.9)
  \psline(2,9.5)(2,7.35)   \psline[linearc=.25](2,7.12)(2,6)(2.5,5.5)(3,6)(3,7.4)
   \psline(2,3.65)(2,0) 
  \psline(3,9.5)(3,7.6)    
\psline[linearc=.25](2,3.9)(2,5)(2.5,5.5)(3,5)(3,3.6)     
\psline(3,3.4)(3,0)
  \psline(4,9.5)(4,7.9)    \psline(4,7.6)(4,3.35)    \psline(4,3.1)(4,1.5)
  \psline(0,0)(0,9.5)      \psline(6,0)(6,9.5)       \psline(7.8,9.5)(7.8,0)
 \psdots(0,9.5)(1,9.5)(2,9.5)(3,9.5)(4,9.5)(5,9.5)(6,9.5)(7.8,9.5)
        (0,0)(1,0)(2,0)(3,0)(4,0)(5,0)(6,0)(7.8,0)
                 \psdots[dotsize=.1](6.65,9.5)(6.95,9.5)(7.25,9.5)
                   \rput(6.95,9.8){$23$'s}
                \psdots[dotsize=.1](6.65,0)(6.95,0)(7.25,0)
                   \rput(6.95,-.3){$23$'s}
                \psdots[dotsize=.1](6.65,4.75)(6.95,4.75)(7.25,4.75)
       \psframe(-.3,9.5)(8.1,10.5)
       \psframe(-.3,0)(8.1,-1)
      \rput(3.65,10.3){$\tau$}
      \rput(3.65,-.8){$\tau'$}
 \rput(0,9.8){$12$}
 \rput(1,9.8){$12$}
 \rput(2,9.8){$15$}  
 \rput(3,9.8){$15$}  
 \rput(4,9.8){$23$}
 \rput(5,9.8){$23$}
 \rput(6,9.8){$23$}
 \rput(7.8,9.8){$23$}
 \rput(0,-.3){$12$}
 \rput(1,-.3){$12$}
 \rput(2,-.3){$14$}
 \rput(3,-.3){$14$}
 \rput(4,-.3){$23$}
 \rput(5,-.3){$23$}
 \rput(6,-.3){$23$}
 \rput(7.8,-.3){$23$}
\end{pspicture}
}}}
\mbox{
\subfigure
{\scalebox{0.5}{
\begin{pspicture}(-.5,-1)(10,10)
  \psline[linearc=.15](5,9.5)(5,5.8)(3,4)(4.9,3)(1,2)(1,1.5)(1,0)
 \psframe[linecolor=white,fillstyle=solid,fillcolor=white](4.5,5.3)(4.7,5.5)
\psline[linearc=.15](4,3.1)(4,1.5)(5,.5)(5,0)
\psframe[linecolor=white,fillstyle=solid,fillcolor=white](4.4,.9)(4.6,1.1)
   \psline[linearc=.15](1,9.5)(1,9)(4.9,8)(3,7)(3.6,6.3)(5,5)(5,1.5)(4,.5)(4,0)
\psframe[linecolor=white,fillstyle=solid,fillcolor=white](1.9,8.65)(2.1,8.85)
\psframe[linecolor=white,fillstyle=solid,fillcolor=white](2.9,8.4)(3.1,8.6)
\psframe[linecolor=white,fillstyle=solid,fillcolor=white](3.9,8.15)(4.1,8.35)
\psframe[linecolor=white,fillstyle=solid,fillcolor=white](1.9,6.7)(2.1,6.9)
\psframe[linecolor=white,fillstyle=solid,fillcolor=white](2.9,6.4)(3.1,6.6)
\psframe[linecolor=white,fillstyle=solid,fillcolor=white](3.85,5.8)(4.15,6.05)
\psframe[linecolor=white,fillstyle=solid,fillcolor=white](1.9,4.2)(2.1,4.4)
\psframe[linecolor=white,fillstyle=solid,fillcolor=white](2.9,4.5)(3.1,4.7)
\psframe[linecolor=white,fillstyle=solid,fillcolor=white](3.9,4.8)(4.15,5.3)
\psframe[linecolor=white,fillstyle=solid,fillcolor=white](1.9,2.1)(2.1,2.4)
\psframe[linecolor=white,fillstyle=solid,fillcolor=white](2.9,2.4)(3.1,2.6)
\psframe[linecolor=white,fillstyle=solid,fillcolor=white](3.9,2.7)(4.1,2.9)
  \psline[linearc=.25](2,9.5)(2,8.2)(2.5,7.7)(3,8.2)(3,9.5)
 \psline[linearc=.25](2,0)(2,2.8)(2.5,3.4)(3,2.8)(3,0)
  \psline(4,9.5)(4,7.65)    \psline(4,7.4)(4,3.6)    \psline(4,3.25)(4,1.5)
  \psline(0,0)(0,9.5)      \psline(6,0)(6,9.5)       \psline(7.8,9.5)(7.8,0)
 \psdots(0,9.5)(1,9.5)(2,9.5)(3,9.5)(4,9.5)(5,9.5)(6,9.5)(7.8,9.5)
        (0,0)(1,0)(2,0)(3,0)(4,0)(5,0)(6,0)(7.8,0)
                 \psdots[dotsize=.1](6.65,9.5)(6.95,9.5)(7.25,9.5)
                   \rput(6.95,9.8){$23$'s}
                \psdots[dotsize=.1](6.65,0)(6.95,0)(7.25,0)
                   \rput(6.95,-.3){$23$'s}
                \psdots[dotsize=.1](6.65,4.75)(6.95,4.75)(7.25,4.75)
       \psframe(-.3,9.5)(8.1,10.5)
       \psframe(-.3,0)(8.1,-1)
      \rput(3.65,10.3){$\tau$}
      \rput(3.65,-.8){$\tau'$}
 \rput(0,9.8){$12$}
 \rput(1,9.8){$12$}
 \rput(2,9.8){$15$}  
 \rput(3,9.8){$15$}  
 \rput(4,9.8){$23$}
 \rput(5,9.8){$23$}
 \rput(6,9.8){$23$}
 \rput(7.8,9.8){$23$}
 \rput(0,-.3){$12$}
 \rput(1,-.3){$12$}
 \rput(2,-.3){$14$}
 \rput(3,-.3){$14$}
 \rput(4,-.3){$23$}
 \rput(5,-.3){$23$}
 \rput(6,-.3){$23$}
 \rput(7.8,-.3){$23$}
\psellipse[linestyle=dotted](4,7)(1,1.4)
\psellipse[linestyle=dotted](4,4)(1,1.4)
\end{pspicture}
}}%
\qquad \qquad
\subfigure
{
\scalebox{0.5}{
\begin{pspicture}(-.5,-1)(10,10)
  \psline[linearc=.15](5,9.5)(5,6)(4,5)(4,1.5)(5,.5)(5,0)
 \psframe[linecolor=white,fillstyle=solid,fillcolor=white](4.4,5.4)(4.6,5.6)
\psline[linearc=.15](4,3.1)(4,1.5)(5,.5)(5,0)
\psframe[linecolor=white,fillstyle=solid,fillcolor=white](4.4,.9)(4.6,1.1)
   \psline[linearc=.15](1,9.5)(1,9)(3.5,8.39)(3.5,2.62)(1,2)(1,1.5)(1,0)
\psframe[linecolor=white,fillstyle=solid,fillcolor=white](1.9,8.65)(2.1,8.85)
\psframe[linecolor=white,fillstyle=solid,fillcolor=white](2.9,8.4)(3.1,8.6)
\psframe[linecolor=white,fillstyle=solid,fillcolor=white](1.9,2.1)(2.1,2.4)
\psframe[linecolor=white,fillstyle=solid,fillcolor=white](2.9,2.4)(3.1,2.6)
  \psline[linearc=.25](2,9.5)(2,8.2)(2.5,7.7)(3,8.2)(3,9.5)
 \psline[linearc=.25](2,0)(2,2.8)(2.5,3.4)(3,2.8)(3,0)
  \psline[linearc=.25](4,9.5)(4,6)(5,5)(5,1.5)(4,.5)(4,0)
  \psline(0,0)(0,9.5)      \psline(6,0)(6,9.5)       \psline(7.8,9.5)(7.8,0)
 \psdots(0,9.5)(1,9.5)(2,9.5)(3,9.5)(4,9.5)(5,9.5)(6,9.5)(7.8,9.5)
        (0,0)(1,0)(2,0)(3,0)(4,0)(5,0)(6,0)(7.8,0)
                 \psdots[dotsize=.1](6.65,9.5)(6.95,9.5)(7.25,9.5)
                   \rput(6.95,9.8){$23$'s}
                \psdots[dotsize=.1](6.65,0)(6.95,0)(7.25,0)
                   \rput(6.95,-.3){$23$'s}
                \psdots[dotsize=.1](6.65,4.75)(6.95,4.75)(7.25,4.75)
       \psframe(-.3,9.5)(8.1,10.5)
       \psframe(-.3,0)(8.1,-1)
      \rput(3.65,10.3){$\tau$}
      \rput(3.65,-.8){$\tau'$}
 \rput(0,9.8){$12$}
 \rput(1,9.8){$12$}
 \rput(2,9.8){$15$}  
 \rput(3,9.8){$15$}  
 \rput(4,9.8){$23$}
 \rput(5,9.8){$23$}
 \rput(6,9.8){$23$}
 \rput(7.8,9.8){$23$}
 \rput(0,-.3){$12$}
 \rput(1,-.3){$12$}
 \rput(2,-.3){$14$}
 \rput(3,-.3){$14$}
 \rput(4,-.3){$23$}
 \rput(5,-.3){$23$}
 \rput(6,-.3){$23$}
 \rput(7.8,-.3){$23$}
\end{pspicture}
}}}
\mbox{
\subfigure
{\scalebox{0.5}{
 \begin{pspicture}(-.5,-1)(10,10)
  \psline(1,9.5)(1,0)      \psline(4,9.5)(4,0)       \psline(5,9.5)(5,0)
   \psline(0,0)(0,9.5)      \psline(6,0)(6,9.5)       \psline(7.8,9.5)(7.8,0)
\psline[linearc=.25](3,9.5)(3,5)(2.5,4.5)(2,5)(2,9.5)
\psline[linearc=.25](3,0)(3,4)(2.5,4.5)(2,4)(2,0)
 \psdots(0,9.5)(1,9.5)(2,9.5)(3,9.5)(4,9.5)(5,9.5)(6,9.5)(7.8,9.5)
        (0,0)(1,0)(2,0)(3,0)(4,0)(5,0)(6,0)(7.8,0)
                 \psdots[dotsize=.1](6.65,9.5)(6.95,9.5)(7.25,9.5)
                   \rput(6.95,9.8){$23$'s}
                \psdots[dotsize=.1](6.65,0)(6.95,0)(7.25,0)
                   \rput(6.95,-.3){$23$'s}
                \psdots[dotsize=.1](6.65,4.75)(6.95,4.75)(7.25,4.75)
       \psframe(-.3,9.5)(8.1,10.5)
       \psframe(-.3,0)(8.1,-1)
      \rput(3.65,10.3){$\tau$}
      \rput(3.65,-.8){$\tau'$}
 \rput(0,9.8){$12$}
 \rput(1,9.8){$12$}
 \rput(2,9.8){$14$}  
 \rput(3,9.8){$14$}  
 \rput(4,9.8){$23$}
 \rput(5,9.8){$23$}
 \rput(6,9.8){$23$}
 \rput(7.8,9.8){$23$}
 \rput(0,-.3){$12$}
 \rput(1,-.3){$12$}
 \rput(2,-.3){$14$}
 \rput(3,-.3){$14$}
 \rput(4,-.3){$23$}
 \rput(5,-.3){$23$}
 \rput(6,-.3){$23$}
 \rput(7.8,-.3){$23$}
\end{pspicture}
}}
\qquad \qquad
\subfigure
{\scalebox{0.5}{
\begin{pspicture}(-.5,-1)(10,10)
  \psline(1,9.5)(1,0)      \psline(4,9.5)(4,0)       \psline(5,9.5)(5,0)
   \psline(0,0)(0,9.5)      \psline(6,0)(6,9.5)       \psline(7.8,9.5)(7.8,0)
    \psline(3,9.5)(3,0)      \psline(2,9.5)(2,0)
 \psdots(0,9.5)(1,9.5)(2,9.5)(3,9.5)(4,9.5)(5,9.5)(6,9.5)(7.8,9.5)
        (0,0)(1,0)(2,0)(3,0)(4,0)(5,0)(6,0)(7.8,0)
                 \psdots[dotsize=.1](6.65,9.5)(6.95,9.5)(7.25,9.5)
                   \rput(6.95,9.8){$23$'s}
                \psdots[dotsize=.1](6.65,0)(6.95,0)(7.25,0)
                   \rput(6.95,-.3){$23$'s}
                \psdots[dotsize=.1](6.65,4.75)(6.95,4.75)(7.25,4.75)
       \psframe(-.3,9.5)(8.1,10.5)
       \psframe(-.3,0)(8.1,-1)
      \rput(3.65,10.3){$\tau$}
      \rput(3.65,-.8){$\tau'$}
 \rput(0,9.8){$12$}
 \rput(1,9.8){$12$}
 \rput(2,9.8){$15$}  
 \rput(3,9.8){$15$}  
 \rput(4,9.8){$23$}
 \rput(5,9.8){$23$}
 \rput(6,9.8){$23$}
 \rput(7.8,9.8){$23$}
 \rput(0,-.3){$12$}
 \rput(1,-.3){$12$}
 \rput(2,-.3){$14$}
 \rput(3,-.3){$14$}
 \rput(4,-.3){$23$}
 \rput(5,-.3){$23$}
 \rput(6,-.3){$23$}
 \rput(7.8,-.3){$23$}
\end{pspicture}
}}}
  \caption{Proof of Theorem~\ref{thm:main} for move V}
  \label{fig:Vproof}
\end{figure}

\begin{figure}[htp]
\begin{center}
  \mbox{
 \subfigure[Begin similarly to Figure~\ref{fig:Iproof}]
{\scalebox{0.6}{
 \begin{pspicture}(.5,0)(7,5)
 \psdots(0,4)(.5,4)(1,4)(1.5,4)(2,4)(2.5,4)(4,4)(4.5,4)
        (0,1)(.5,1)(1,1)(1.5,1)(2,1)(2.5,1)(4,1)(4.5,1)
   \psline[linearc=.25](1,4)(1,3)(1.25,2.5)(1.5,3)(1.5,4)
   \psline[linearc=.25](1,1)(1,2)(1.25,2.5)(1.5,2)(1.5,1)
   \psline(2,4)(2,1)
   \psline(2.5,4)(2.5,1)
   \psline(4,4)(4,1)
   \psline(4.5,4)(4.5,1) 
 \psline[linearc=.25](.5,4)(.5,3)(0,2)(0,1)
\psframe[linecolor=white,fillstyle=solid,fillcolor=white](.2,2.4)(.39,2.6)
\psline[linearc=.255](0,4)(0,3)(.5,2)(.5,1)
           \psframe(-.3,1)(4.8,0)
           \psframe(-.3,4)(4.8,5)
    \rput(0,4.3){$12$}
 \rput(0.5,4.3){$12$}
 \rput(1,4.3){$15$}
 \rput(1.5,4.3){$15$}
 \rput(2,4.3){$23$}
 \rput(2.5,4.3){$23$}
 \rput(4,4.3){$23$}
 \rput(4.5,4.3){$23$}
                   \rput(3.2,4.3){$23$'s}
 \rput(0,.7){$12$}
 \rput(0.5,.7){$12$}
 \rput(1,.7){$14$}
 \rput(1.5,.7){$14$}
 \rput(2,.7){$23$}
 \rput(2.5,.7){$23$}
 \rput(4,.7){$23$}
 \rput(4.5,.7){$23$}
                   \rput(3.2,.7){$23$'s}
  \psdots[dotsize=0.1](3,4)(3.3,4)(3.6,4)(3,1)(3.3,1)(3.6,1)
                      (3,2.5)(3.3,2.5)(3.6,2.5)
     \rput(2.5,4.8){$\tau$}
     \rput(2.5,.2){$\tau'$} 
\end{pspicture}
}
}
\qquad
\qquad
\qquad
\subfigure[then use isotopy]{\scalebox{0.5}{
\begin{pspicture}(-.5,-2)(6,5)
 \psdots(0,4)(.5,4)(2,4)(2.5,4)(3,4)(3.5,4)(5,4)(5.5,4)
        (0,-1)(.5,-1)(2,-1)(2.5,-1)(3,-1)(3.5,-1)(5,-1)(5.5,-1)
   \psline[linearc=.25](2,4)(2,3)(2.25,2.5)(2.5,3)(2.5,4)
   \psline(3,4)(3,-1)
   \psline(3.5,4)(3.5,-1)
   \psline(5,4)(5,-1)
   \psline(5.5,4)(5.5,-1) 
\psline[linearc=.25](.5,4)(.5,3.6)(1.4,3.6)(1.4,2.5)(1.25,2)(.9,2.5)(.9,3.2)(.5,3.2)(.5,2)(0,1)(0,-1)
\psframe[linecolor=white,fillstyle=solid,fillcolor=white](.2,1.4)(.39,1.6)
\psline[linearc=.25](0,4)(0,2)(.5,1)(.5,-1)
\psframe[linecolor=white,fillstyle=solid,fillcolor=white](.4,-.7)(.6,.9)
\psline[linearc=.15](2,-1)(2,-.4)(.5,-.4)(.05,.2)(.5,.8)(.8,.8)(.8,1.8)(1.1,2.3)(1.4,1.8)(1.4,.4)(1.1,.4)(-1.5,.2)(1.1,0)(2.5,0)(2.5,-1)
\psframe[linecolor=white,fillstyle=solid,fillcolor=white](.4,-.7)(.6,.2)
\psline(.5,.2)(.5,-1)  \psline(.5,.65)(.5,.45)
           \psframe(-.3,-1)(5.8,-2)
           \psframe(-.3,4)(5.8,5)
    \rput(0,4.3){$12$}
 \rput(0.5,4.3){$12$}
 \rput(2,4.3){$15$}
 \rput(2.5,4.3){$15$}
 \rput(3,4.3){$23$}
 \rput(3.5,4.3){$23$}
 \rput(5,4.3){$23$}
 \rput(5.5,4.3){$23$}
                   \rput(4.2,4.3){$23$'s}
 \rput(0,-1.3){$12$}
 \rput(0.5,-1.3){$12$}
 \rput(2,-1.3){$14$}
 \rput(2.5,-1.3){$14$}
 \rput(3,-1.3){$23$}
 \rput(3.5,-1.3){$23$}
 \rput(5,-1.3){$23$}
 \rput(5.5,-1.3){$23$}
                   \rput(4.2,-1.3){$23$'s}
  \psdots[dotsize=0.1](4,4)(4.3,4)(4.6,4)(4,-1)(4.3,-1)(4.6,-1)
                      (4,1.5)(4.3,1.5)(4.6,1.5)
\rput(1.7,2){$24$}
     \rput(3,4.8){$\tau$}
     \rput(3,-1.8){$\tau'$} 
\end{pspicture}
}
}
}
 \mbox{
 \subfigure[then apply reverse circumcision]{\scalebox{0.5}{
\begin{pspicture}(-.5,-2)(8,5)
 \psdots(0,4)(.5,4)(2,4)(2.5,4)(3,4)(3.5,4)(5,4)(5.5,4)
        (0,-1)(.5,-1)(2,-1)(2.5,-1)(3,-1)(3.5,-1)(5,-1)(5.5,-1)
   \psline[linearc=.25](2,4)(2,3)(2.25,2.5)(2.5,3)(2.5,4)
   \psline(3,4)(3,-1)
   \psline(3.5,4)(3.5,-1)
   \psline(5,4)(5,-1)
   \psline(5.5,4)(5.5,-1) 
\psline[linearc=.15](.5,4)(.5,3.6)(1.4,3.6)
(1.4,.5)(1,.5)(0,.1)(1,-.3)(2.5,-.3)(2.5,-1)
\psline[linearc=.15](2,-1)(2,-.6)(1,-.6)(.5,-.4)
(.05,.1)(.5,.7)(1,.9)(1,3.2)(.5,3.2)(.5,2)(0,1)(0,-1)
\psframe[linecolor=white,fillstyle=solid,fillcolor=white](.2,1.4)(.4,1.6)
\psline[linearc=.25](0,4)(0,2)(.5,1)(.5,.8)
\psframe[linecolor=white,fillstyle=solid,fillcolor=white](.4,0)(.6,-.8)
\psline(.5,1)(.5,.8)  \psline(.5,.55)(.5,.41) \psline(.5,.2)(.5,-1)
           \psframe(-.3,-1)(5.8,-2)
           \psframe(-.3,4)(5.8,5)
    \rput(0,4.3){$12$}
 \rput(0.5,4.3){$12$}
 \rput(2,4.3){$15$}
 \rput(2.5,4.3){$15$}
 \rput(3,4.3){$23$}
 \rput(3.5,4.3){$23$}
 \rput(5,4.3){$23$}
 \rput(5.5,4.3){$23$}
                   \rput(4.2,4.3){$23$'s}
 \rput(0,-1.3){$12$}
 \rput(0.5,-1.3){$12$}
 \rput(2,-1.3){$14$}
 \rput(2.5,-1.3){$14$}
 \rput(3,-1.3){$23$}
 \rput(3.5,-1.3){$23$}
 \rput(5,-1.3){$23$}
 \rput(5.5,-1.3){$23$}
                   \rput(4.2,-1.3){$23$'s}
  \psdots[dotsize=0.1](4,4)(4.3,4)(4.6,4)(4,-1)(4.3,-1)(4.6,-1)
                      (4,1.5)(4.3,1.5)(4.6,1.5)
\rput(1.7,.6){$24$}
     \rput(3,4.8){$\tau$}
     \rput(3,-1.8){$\tau'$}
 \rput(1.9,2){$14$}
\psframe[linecolor=white,fillstyle=solid,fillcolor=white](.8,1.7)(1.6,2)
\psellipse(1.2,2)(.4,.3) 
\psframe[linecolor=white,fillstyle=solid,fillcolor=white](.9,2)(1.1,2.4)
\psframe[linecolor=white,fillstyle=solid,fillcolor=white](1.3,2)(1.5,2.4)
\psline(1,2.5)(1,1.9)   \psline(1.4,2.5)(1.4,1.9)
\end{pspicture}
}
}
\qquad
\qquad
\qquad
 \subfigure[apply a sequence of $\CP$ moves]{\scalebox{0.5}{
\begin{pspicture}(-.5,-2)(6,5)
 \psdots(0,4)(.5,4)(2,4)(2.5,4)(3,4)(3.5,4)(5,4)(5.5,4)
        (0,-1)(.5,-1)(2,-1)(2.5,-1)(3,-1)(3.5,-1)(5,-1)(5.5,-1)
   \psline[linearc=.25](2,4)(2,3)(2.25,2.5)(2.5,3)(2.5,4)
   \psline(3,4)(3,-1)
   \psline(3.5,4)(3.5,-1)
   \psline(5,4)(5,-1)
   \psline(5.5,4)(5.5,-1) 
\psline[linearc=.15](.5,4)(.5,3.6)(1.4,3.6)(1.4,.5)(1,.5)
(0,.1)(1,-.3)(2.5,-.3)(2.5,-1)
\psline[linearc=.15](2,-1)(2,-.6)(1,-.6)(.5,-.4)(.05,.1)(.5,.7)
(1,.9)(1,3.2)(.5,3.2)(.5,2)(0,1)(0,-1)
\psframe[linecolor=white,fillstyle=solid,fillcolor=white](.2,1.4)(.4,1.6)
\psline[linearc=.25](0,4)(0,2)(.5,1)(.5,.8)
\psframe[linecolor=white,fillstyle=solid,fillcolor=white](.4,0)(.6,-.8)
\psline(.5,1)(.5,.8)  \psline(.5,.55)(.5,.41) \psline(.5,.2)(.5,-1)
           \psframe(-.3,-1)(5.8,-2)
           \psframe(-.3,4)(5.8,5)
    \rput(0,4.3){$12$}
 \rput(0.5,4.3){$12$}
 \rput(2,4.3){$15$}
 \rput(2.5,4.3){$15$}
 \rput(3,4.3){$23$}
 \rput(3.5,4.3){$23$}
 \rput(5,4.3){$23$}
 \rput(5.5,4.3){$23$}
                   \rput(4.2,4.3){$23$'s}
 \rput(0,-1.3){$12$}
 \rput(0.5,-1.3){$12$}
 \rput(2,-1.3){$14$}
 \rput(2.5,-1.3){$14$}
 \rput(3,-1.3){$23$}
 \rput(3.5,-1.3){$23$}
 \rput(5,-1.3){$23$}
 \rput(5.5,-1.3){$23$}
                   \rput(4.2,-1.3){$23$'s}
  \psdots[dotsize=0.1](4,4)(4.3,4)(4.6,4)(4,-1)(4.3,-1)(4.6,-1)
                      (4,1.5)(4.3,1.5)(4.6,1.5)
\rput(1.7,.6){$24$}
     \rput(3,4.8){$\tau$}
     \rput(3,-1.8){$\tau'$}
 \rput(2.6,1.5){$14$}
\psframe[linecolor=white,fillstyle=solid,fillcolor=white](.8,1.65)(5.6,2)
\psellipse(3.25,2.)(2.6,.3) 
\psframe[linecolor=white,fillstyle=solid,fillcolor=white](.9,2)(1.1,2.4)
\psframe[linecolor=white,fillstyle=solid,fillcolor=white](1.3,2)(1.5,2.4)
\psframe[linecolor=white,fillstyle=solid,fillcolor=white](2.9,2)(3.1,2.4)
\psframe[linecolor=white,fillstyle=solid,fillcolor=white](3.4,2)(3.6,2.4)
\psframe[linecolor=white,fillstyle=solid,fillcolor=white](4.9,2)(5.1,2.4)
\psframe[linecolor=white,fillstyle=solid,fillcolor=white](5.4,2)(5.6,2.4)
\psline(1,2.5)(1,1.96)   \psline(1.4,2.5)(1.4,1.91)
\psline(3,2.5)(3,1.89)   \psline(3.5,2.5)(3.5,1.9)
\psline(5,2.5)(5,1.91)   \psline(5.5,2.5)(5.5,1.96)
\end{pspicture}
}
}
}
\mbox{
 \subfigure[flip the circumcising  sheet over the diagram]{\scalebox{0.5}{
\begin{pspicture}(-.5,-2)(8,5)
 \psdots(0,4)(.5,4)(2,4)(2.5,4)(3,4)(3.5,4)(5,4)(5.5,4)
        (0,-1)(.5,-1)(2,-1)(2.5,-1)(3,-1)(3.5,-1)(5,-1)(5.5,-1)
   \psline[linearc=.25](2,4)(2,3)(2.25,2.5)(2.5,3)(2.5,4)
   \psline(3,4)(3,-1)
   \psline(3.5,4)(3.5,-1)
   \psline(5,4)(5,-1)
   \psline(5.5,4)(5.5,-1) 
\psline[linearc=.15](.5,4)(.5,3.6)(1.4,3.6)(1.4,.5)(1,.5)
(0,.1)(1,-.3)(2.5,-.3)(2.5,-1)
\psline[linearc=.15](2,-1)(2,-.6)(1,-.6)(.5,-.4)(.05,.1)
(.5,.7)(1,.9)(1,3.2)(.5,3.2)(.5,2)(0,1)(0,-1)
\psframe[linecolor=white,fillstyle=solid,fillcolor=white](.2,1.4)(.4,1.6)
\psline[linearc=.25](0,4)(0,2)(.5,1)(.5,.8)
\psframe[linecolor=white,fillstyle=solid,fillcolor=white](.4,0)(.6,-.8)
\psline(.5,1)(.5,.8)  \psline(.5,.55)(.5,.41) \psline(.5,.2)(.5,-1)
           \psframe(-.3,-1)(5.8,-2)
           \psframe(-.3,4)(5.8,5)
    \rput(0,4.3){$24$}
 \rput(0.5,4.3){$24$}
 \rput(2,4.3){$45$}
 \rput(2.5,4.3){$45$}
 \rput(3,4.3){$23$}
 \rput(3.5,4.3){$23$}
 \rput(5,4.3){$23$}
 \rput(5.5,4.3){$23$}
                   \rput(4.2,4.3){$23$'s}
 \rput(0,-1.3){$12$}
 \rput(0.5,-1.3){$12$}
 \rput(2,-1.3){$14$}
 \rput(2.5,-1.3){$14$}
 \rput(3,-1.3){$23$}
 \rput(3.5,-1.3){$23$}
 \rput(5,-1.3){$23$}
 \rput(5.5,-1.3){$23$}
                   \rput(4.2,-1.3){$23$'s}
  \psdots[dotsize=0.1](4,4)(4.3,4)(4.6,4)(4,-1)(4.3,-1)(4.6,-1)
                      (4,1.5)(4.3,1.5)(4.6,1.5)
\rput(1.7,.6){$24$}
     \rput(3,4.8){$\tau$}
     \rput(3,-1.8){$\tau'$}
 \rput(-.5,2.3){$14$}
\psframe[linecolor=white,fillstyle=solid,fillcolor=white](-.1,2)(.6,2.2)
\psellipse(.25,2.3)(.5,.25) 
\psframe[linecolor=white,fillstyle=solid,fillcolor=white](-.1,2.2)(.1,2.7)
\psframe[linecolor=white,fillstyle=solid,fillcolor=white](.4,2.2)(.6,2.7)
\psline(0,3)(0,2.2)   \psline(.5,3)(.5,2.2)
\end{pspicture}
}
}
\qquad
\qquad
\qquad
\subfigure[and circumcise. Now finish the proof similarly to move I.]{\scalebox{0.5}{
\begin{pspicture}(-.5,-2)(6,5)
 \psdots(0,4)(.5,4)(2,4)(2.5,4)(3,4)(3.5,4)(5,4)(5.5,4)
        (0,-1)(.5,-1)(2,-1)(2.5,-1)(3,-1)(3.5,-1)(5,-1)(5.5,-1)
   \psline[linearc=.25](2,4)(2,3)(2.25,2.5)(2.5,3)(2.5,4)
   \psline(3,4)(3,-1)
   \psline(3.5,4)(3.5,-1)
   \psline(5,4)(5,-1)
   \psline(5.5,4)(5.5,-1) 
\psline[linearc=.15](.5,4)(.5,3.6)(1.4,3.6)(1.4,.5)
(1,.5)(0,.1)(1,-.3)(2.5,-.3)(2.5,-1)
\psline[linearc=.15](2,-1)(2,-.6)(1,-.6)(.5,-.4)(.05,.1)(.5,.7)(1,.9)
(1,3.2)(.5,3.2)(.5,3)(.25,2.5)(0,3)(0,4)
\psline[linearc=.25](0,2)(.25,2.8)(.5,2)(0,1)(0,-1)
\psframe[linecolor=white,fillstyle=solid,fillcolor=white](.2,1.4)(.4,1.6)
\psline[linearc=.25](.5,2)(.25,2.8)(0,2)(.5,1)(.5,.8)
\psframe[linecolor=white,fillstyle=solid,fillcolor=white](.4,0)(.6,-.8)
\psline(.5,1)(.5,.8)  \psline(.5,.55)(.5,.41) \psline(.5,.2)(.5,-1)
           \psframe(-.3,-1)(5.8,-2)
           \psframe(-.3,4)(5.8,5)
    \rput(0,4.3){$24$}
 \rput(0.5,4.3){$24$}
 \rput(2,4.3){$45$}
 \rput(2.5,4.3){$45$}
 \rput(3,4.3){$23$}
 \rput(3.5,4.3){$23$}
 \rput(5,4.3){$23$}
 \rput(5.5,4.3){$23$}
                   \rput(4.2,4.3){$23$'s}
 \rput(0,-1.3){$12$}
 \rput(0.5,-1.3){$12$}
 \rput(2,-1.3){$14$}
 \rput(2.5,-1.3){$14$}
 \rput(3,-1.3){$23$}
 \rput(3.5,-1.3){$23$}
 \rput(5,-1.3){$23$}
 \rput(5.5,-1.3){$23$}
                   \rput(4.2,-1.3){$23$'s}
  \psdots[dotsize=0.1](4,4)(4.3,4)(4.6,4)(4,-1)(4.3,-1)(4.6,-1)
                      (4,1.5)(4.3,1.5)(4.6,1.5)
\rput(1.7,.6){$24$}
     \rput(3,4.8){$\tau$}
     \rput(3,-1.8){$\tau'$}
\end{pspicture}
}
}
}
\caption{Proof of Theorem~\ref{thm:main} for move II}
\label{fig:IIproof}
\end{center}
\end{figure}

\np

\Addresses\recd
\end{document}